\documentclass[10pt, reqno]{amsart}
\usepackage{hyperref}
\hypersetup{
	colorlinks=true,
	citecolor=blue,
	linkcolor=blue,
	filecolor=magenta,      
	urlcolor=cyan,
}

\usepackage[margin=0.5in]{geometry}
\usepackage{setspace}
\usepackage[capitalize,noabbrev,nameinlink]{cleveref}
\usepackage{fullpage}
\usepackage{amsmath, amsthm, amssymb, graphicx, dsfont, stmaryrd, extarrows, enumitem}
\usepackage[T1]{fontenc}
\usepackage{textcomp}
\usepackage{lmodern}
\usepackage{microtype} 
\usepackage{amscd}
\usepackage{mathrsfs}
\usepackage{tikz-cd}
\usetikzlibrary{matrix,patterns}
\usepackage[colorinlistoftodos]{todonotes}
\usepackage{color}
\usepackage{bbm}
\usepackage{verbatim}
\usepackage[mathscr]{euscript}
\usepackage{bm}

\numberwithin{equation}{section}
\theoremstyle{plain}
\newtheorem{thm}[equation]{Theorem}
\newtheorem*{thm*}{Theorem}

\newtheorem{prop}[equation]{Proposition}
    
\newtheorem{cor}[equation]{Corollary}       
\newtheorem{lem}[equation]{Lemma}

\theoremstyle{definition}

\newtheorem{ex}[equation]{Example}

\newtheorem{rem}[equation]{Remark}

\newtheorem{chunk}[equation]{}

\newcommand{\Z}{\mathbb{Z}}

\newcommand{\mf}{\mathfrak}
\newcommand{\Hom}{\mathrm{Hom}}

\newcommand{\msf}[1]{\mathsf{#1}}

\newcommand{\mc}[1]{\mathcal{#1}}
\newcommand{\mrm}[1]{\mathrm{#1}}
\newcommand{\mbb}[1]{\mathbb{#1}}
\newcommand{\scr}[1]{\mathscr{#1}}

\newcommand{\T}{\mathsf{T}}
\newcommand{\U}{\mathsf{U}}

\newcommand{\1}{\mathds{1}}

\renewcommand{\mod}[1]{\mathrm{Mod}_{#1}}

\newcommand{\A}{\mathsf{A}}
\newcommand{\B}{\mathsf{B}}

\newcommand{\D}{\mathsf{D}}

\renewcommand{\mod}[1]{\msf{mod}(#1)}
\newcommand{\Mod}[1]{\msf{Mod}(#1)}
\newcommand{\Flat}[1]{\msf{Flat}(#1)}
\newcommand{\ab}{\msf{Ab}}
\newcommand{\ti}{\textit}
\renewcommand{\tilde}[1]{\widetilde{#1}}
\newcommand{\rlim}{\varinjlim}
\newcommand{\llim}{\varprojlim}
\renewcommand{\t}{\text}
\renewcommand{\c}{\mrm{c}}

\newcommand{\op}{\mrm{op}}

\renewcommand{\Flat}[1]{\msf{Flat}(#1)}

\newcommand{\y}{\msf{y}}
\newcommand{\stgf}[1]{\underline{\smash{\msf{GFlatCot}}}(#1)}
\newcommand{\gfc}[1]{\msf{GFlatCot}(#1)}

\newcommand{\stgp}[1]{\underline{\smash{\msf{Gproj}}}(#1)}
\newcommand{\stgP}[1]{\underline{\smash{\msf{GProj}}}(#1)}
\newcommand{\ul}[1]{\underline{\smash{#1}}}

\title{Purity, ascent, and periodicity for Gorenstein flat cotorsion modules}

\begin{document}
\author{Isaac Bird}
\address{Department of Algebra, Faculty of Mathematics and Physics, Charles University in Prague, Sokolovsk\'{a} 83, 186 75 Praha, Czech Republic}
\email{bird@karlin.mff.cuni.cz}

\maketitle

\begin{abstract}
We investigate purity within the Frobenius category of Gorenstein flat cotorsion modules, which can be seen as an infinitely generated analogue of the Frobenius category of Gorenstein projective objects. As such, the associated stable category can be viewed as an alternative approach to a big singularity category, which is equivalent to Krause's when the ring is Gorenstein.

We study the pure structure of the stable category, and show it is fundamentally related to the pure structure of the Gorenstein flat modules. Following that, we give conditions for extension of scalars to preserve Gorenstein flat cotorsion modules. In this case, one obtains an induced triangulated functor on the stable categories. We show that under mild conditions that these functors preserve the pure structure, both on the triangulated and module category level.

Along the way, we consider particular phenomena over commutative rings, the cumulation of which is an extension of Kn\"{o}rrer periodicity, giving a triangulated equivalence between Krause's big singularity categories for a complete hypersurface singularity and its twofold double-branched cover.
\end{abstract}

\setcounter{tocdepth}{1}
\tableofcontents

\section{Introduction}
Understanding the global dimension of a ring, or the smoothness of a scheme, is a fundamental question in algebra and geometry. One of the more celebrated constructions which measures this regularity is the singularity category, as introduced by Buchweitz \cite{buchweitz}, and independently later by Orlov \cite{orlov2}. For a ring $R$, its singularity category is defined to be the Verdier quotient $\msf{D}_{\t{sg}}(R):=\msf{D}^{\t{b}}(\mod{R})/\msf{perf}(R)$, and it is clear that, whenever $R$ is sufficiently nice, this vanishes if and only if $R$ has finite global dimension.

However, this is far from being the only triangulated category which measures regularity. If one approaches from the direction of Gorenstein homological algebra, one can consider the stable category $\stgp{R}$ of finitely presented Gorenstein projective modules. Again, when $R$ is sufficiently nice,  this category vanishes if and only if $R$ has finite global dimension. It is possible to compare $\stgp{R}$ and $\msf{D}_{\t{sg}}(R)$. Over any ring there is an embedding $\stgp{R}\to\msf{D}_{\t{sg}}(R)$, and one of the striking results of \cite{buchweitz} is that this functor is an equivalence of categories if, and it transpires only if, $R$ is an Iwanaga-Gorenstein ring.

Both $\msf{D}_{\t{sg}}(R)$ and $\stgp{R}$ are essentially small triangulated categories, and firmly live in the world of bounded and finite objects. While for many problems this finite world may be all that one needs, for certain questions, such as those considered in this paper, they do not suffice. Instead, one wishes to consider an unbounded analogue. 

In this paper, we investigate the triangulated category which provides the unbounded analogue to $\stgp{R}$ when one approaches the problem from the perspective of purity. The category which we take for this unbounded analogue is $\stgf{R}$, the stable category of Gorenstein flat and cotorsion modules, or equivalently $\msf{K}_{\t{tac}}(\msf{FlatCot}(R))$, the homotopy category of totally acyclic complexes of flat and cotorsion modules. We study the pure structure of $\stgf{R}$, and how it behaves across changes of rings. Before explicitly defining $\stgf{R}$, giving a justification as to why it is our object of study, and stating the main results, we motivate the pure structure of a triangulated category. 

The pure structure on a compactly generated triangulated category contains a substantial amount of information, which is of interest across several areas of algebra and topology. For example, it completely determines the smashing subcategories of the category for stable homotopy theory, see \cite{krsmash, krcoh}, encodes information about (co)silting objects for representation theory, see \cite{AMV, laking, lv, MarksVitoria}, determines the rank functions on a triangulated category, see \cite{conde2022functorial}, and, in the case of tensor triangulated categories, information about the Balmer spectrum, see \cite{bks,definablefunctors, birdwilliamsonhomological}. Thus understanding the pure structure really reveals a substantial amount of desirable information. 

This pure structure consists of two parts. Firstly, there are the pure injective objects. In general, it is impossible to classify indecomposable non-finitely generated, or non-compact, objects, as there is usually a proper class of them. However, there is always only a set of indecomposable pure injective objects, up to isomorphism. The second part of the pure structure enables one to recover the category, up to purity, from the indecomposable pure injectives. This second part consists of definable subcategories, which can be thought of as vanishing sets of particularly well behaved functors. Together, one can form a topological space which contains, and determines, the entirety of the pure structure of the category, called the Ziegler spectrum. Its points are the indecomposable pure injectives, and the definable subcategories parametrise its closed sets. These concepts are equally well defined in module categories, see \cite{krspec,psl} for surveys.

It is through considering the pure structure on modules that one is led to the triangulated category $\stgf{R}$. Assuming our ring is sufficiently nice, if one wants to consider the part of the pure structure that is determined by $\msf{proj}(R)$, one ends up with the definable category $\msf{Flat}(R)$. Similarly one considers the definable category $\msf{GFlat}(R)$ to witness $\msf{Gproj}(R)$, and there is a natural inclusion $\msf{Flat}(R)\subseteq\msf{GFlat}(R)$ which is the `purity' analogue of the inclusion $\msf{proj}(R)\subseteq\msf{Gproj}(R)$. Any pure injective module $X$ is \ti{cotorsion}, that is satisfies $\t{Ext}_{R}^{1}(F,X)=0$ for any flat module $F$. Therefore the pure structure on $\msf{GFlat}(R)$ is completely contained within the category $\gfc{R}$ consisting of modules which are simultaneously Gorenstein flat and cotorsion; likewise the category $\msf{FlatCot}(R)$ consisting of modules which are simultaneously flat and cotorsion completely determines the pure structure of $\msf{Flat}(R)$. Yet $\gfc{R}$ happens to be a Frobenius exact category, whose injective objects are $\msf{FlatCot}(R)$, and thus solely by considering the pure structure induced by the inclusion $\msf{proj}(R)\subseteq\msf{Gproj}(R)$ we are led to $\stgf{R}$.

At this point it would be amiss not to mention the already well established big analogue of $\msf{D}_{\t{sg}}(R)$. Originally introduced by Krause \cite{krstable}, and generalised in \cite{becker, spure}, if one uses the injective $R$-modules $\msf{Inj}(R)$ over particularly nice rings, one obtains a recollement
\[
\begin{tikzcd}
\msf{K}_{\t{ac}}(\msf{Inj}(R)) \arrow[r] \arrow[r, leftarrow, shift left=1ex] \arrow[r, leftarrow, shift right=1ex] & \msf{K}(\msf{Inj}(R)) \arrow[r] \arrow[r, leftarrow, shift left=1ex] \arrow[r, leftarrow, shift right=1ex] & \D(R)
\end{tikzcd}
\]
which, when taking compact objects, yields the localisation sequence
\[
\begin{tikzcd}
\msf{D}_{\t{sg}}(R) \arrow[r, leftarrow] & \msf{D}^{\t{b}}(\mod{R}) \arrow[r, leftarrow] & \msf{perf}(R)
\end{tikzcd}
\]
giving Buchweitz's original construction, and realising $\msf{K}_{\t{ac}}(\msf{Inj}(R))$ as a `big' singularity category. One would hope that $\stgf{R}$ would behave in relation to $\msf{K}_{\t{ac}}(R)$ in an analogous way that $\stgp{R}$ does in relation to $\msf{D}_{\t{sg}}(R)$, and this is indeed the case: over any ring there is an embedding $\stgf{R}\to\msf{K}_{\t{ac}}(\msf{Inj}(R))$, and whenever $R$ is Gorenstein, the categories $\stgf{R}$ and $\msf{K}_{\t{ac}}(\msf{Inj}(R))$ are actually triangle equivalent.

Let us now discuss the main results of the paper. As mentioned above, from considering the pure structure on the definable classes $\msf{Flat}(R)$ and $\msf{GFlat}(R)$, we get to $\stgf{R}$. It is therefore a natural question as to how the pure structure on $\msf{GFlat}(R)$ and $\stgf{R}$ are related; this question is considered in \cref{purestructureofstable}. For the pure injective objects, this is straightforward: it is shown in \cref{piinst} that an object in $\stgf{R}$ is pure injective if and only if it is pure injective in $\msf{GFlat}(R)$. 

Relating the remainder of the pure structure is slightly more difficult. Firstly, without any assumptions on the ring there is no reason to suppose that $\stgf{R}$ is compactly generated. However, this is also the case for $\msf{K}_{\t{ac}}(\msf{Inj}(R))$, so we only consider the cases when compact generation holds, which includes the nicest cases of $R$ being Iwanaga-Gorenstein or, more generally, being a Noetherian ring admitting a dualising complex. However, identifying the compact objects of $\stgf{R}$ is not straightforward, since, as shown in \cref{coproductsperfect}, the coproducts in $\stgf{R}$ are not those of $\Mod{R}$ unless $R$ is a left perfect ring. Nevertheless, with our assumption of compact generation, the pure triangles of $\stgf{R}$ are related to the pure exact sequences in $\gfc{R}$ in \cref{pureexactrelationship}.

Under hypotheses which include the common cases of $R$ being a complete commutative Gorenstein local ring, or a Gorenstein Artin algebra, we are able to completely determine the pure triangles and definable subcategories of $\stgf{R}$. This leads to the following theorem.

\begin{thm*}[\ref{extensionpureiff}, \ref{zieglerdisjoint}]
Let $R$ be a ring that is cotorsion over itself such that $\msf{Gproj}(R)\subseteq\msf{GFlat}(R)$, there is an equality $\msf{GFlat}(R)=\rlim\msf{Gproj}(R)$ and the objects of $\stgp{R}$ are compact in $\stgf{R}$. Then $\stgf{R}$ is compactly generated and there is a homeomorphism
\[
\msf{Zg}(\msf{GFlat}(R))\simeq \msf{Zg}(\stgf{R})\coprod\msf{Zg}(\Flat{R}).
\]
\end{thm*}

Once one has established a relationship between purity in $\msf{GFlat}(R)$ and $\stgf{R}$, it is natural to wonder how this behaves, for example, with respect to change of rings. In \cref{ascentsection}, we investigate how extension of scalars along a ring homomorphism $R\to S$ transfers the pure structure between $\msf{GFlat}(R)$ and $\msf{GFlat}(S)$, as well as the corresponding question for the categories $\stgf{R}$ and $\stgf{S}$. In order to preserve the pure structure, one requires $S$ to be finitely presented over $R$, as this means $S\otimes_{R}-$ preserves direct limits and products, and thus preserves pure exact sequences and pure injective modules. Yet we also wish to preserve Gorenstein flat modules, which requires $S$ to have finite flat dimension. It transpires that, with these two assumptions, we actually preserve substantially more.

\begin{thm*}[{\ref{triangulatedfunctor}}]
Let $R\to S$ be a morphism of right coherent rings such that $S$ is finitely presented of finite flat dimension as a right $R$-module. Then extension of scalars induces a functor
\[
S\otimes_{R}-\colon\msf{Ch}_{\mrm{tac}}(\msf{FlatCot}(R))\to\msf{Ch}_{\mrm{tac}}(\msf{FlatCot}(S))
\]
between the categories of totally acyclic complexes of flat-cotorsion modules. In particular, it yields a functor $\gfc{R}\to\gfc{S}$ of Frobenius exact categories, and the induced functor
\[
S\otimes_{R}-\colon\stgf{R}\to\stgf{S}
\] 
is triangulated.
\end{thm*}
Following this, further properties of this induced triangulated functor are considered, for instance, in \cref{preservesprodsandcoprods} it is shown that it always preserves products, while under certain non-restrictive circumstances it also preserves coproducts. In the case that $\stgf{R}$ and $\stgf{S}$ are compactly generated, this means that $S\otimes_{R}-$ is a definable functor in the sense of \cite{definablefunctors}. 

As the pure structure, in the sense of pure injective objects and definable subcategories, is preserved, one may wonder whether $S\otimes_{R}-$ induces a map of topological spaces $\msf{Zg}(\msf{GFlat}(R))\to\msf{Zg}(\msf{GFlat}(R))$. One does not get this for free - it requires the functor to be full on pure injective objects. Bearing in mind that the pure injective objects in $\msf{GFlat}(R)$ are the pure injective objects in $\stgf{R}$, we then show that for many pleasant ring homomorphisms it is actually possible to check this fullness assumption on the level of triangulated categories, where there are fewer morphisms.

\begin{thm*}[{\ref{fullonpis}, \ref{conditionsholdforcommutative}}]
Let $R\to S$ be a surjective map of right coherent rings such that $S$ is finitely presented of finite flat dimension over $R$. If the functor $S\otimes_{R}-\colon\stgf{R}\to\stgf{S}$ is full on pure injective objects, then the functor $S\otimes_{R}-\colon\msf{GFlat}(R)\to\msf{GFlat}(S)$ is full on pure injectives, and thus induces a closed and continuous map $\msf{Zg}(\msf{GFlat}(R))\to\msf{Zg}(\msf{GFlat}(S))$.
\end{thm*}

Following this, we further investigate the interplay of extension of scalars on the level of definable, Frobenius, and triangulated categories. Of interest is the kernel of all these functors, which has significance in understanding the Ziegler spectrum of the image. In greatest generality, it is shown in \cref{kernelfrobenius} that under certain assumptions the kernel of $S\otimes_{R}-\colon\gfc{R}\to\gfc{S}$ is itself a Frobenius category whose stable category is colocalising in the kernel of $\stgf{R}\to\stgf{S}$; in many cases these two triangulated categories have the same pure injective objects.

In \cref{commrings}, we study phenomena occurring when $R$ is a commutative Noetherian ring. The first main result shows that in a familiar setting the assumptions of the aforementioned proposition hold, and the kernel of the extension of scalars functor is a Frobenius category.

\begin{thm*}[{\ref{comkerfrob}}]
Let $R$ be a commutative Noetherian ring and let $I$ be an ideal generated by an $R$-sequence. Set $\msf{K}=\{X\in\gfc{R}:R/I\otimes_{R}X=0\}$. Then $\msf{K}$ is closed under flat covers and Gorenstein cotorsion (pre)envelopes. In particular, $\msf{K}$ is a Frobenius category such that $\ul{\msf{K}}=\{Z\in\stgf{R}:R/I\otimes_{R}Z=0\}$.
\end{thm*}

Following this, we investigate Gorenstein flat cotorsion modules over hypersurface singularities. One of the most celebrated results concerning modules over hypersurface singularities is Kn\"{o}rrer periodicity, which states that there is a triangulated equivalence of categories $\msf{D}_{\t{sg}}(R)\simeq \msf{D}_{\t{sg}}(R^{\sharp\sharp})$, where $R^{\sharp\sharp}$ is the two-fold double branched cover of $R$, see \cite{knorrer}. We investigate this phenomenon for the category $\stgf{R}$ which, as $R$ is Gorenstein, is equivalent to Krause's big singularity category $\msf{K}_{\t{ac}}(\msf{Inj}(R))$; it transpires that one obtains the following periodicity result, extending Kn\"{o}rrer's.

\begin{thm*}[{\ref{bigknorrer}}]
Let $S=k[[x_{1},\cdots,x_{n+1}]]$ be a complete hypersurface ring with $k$ algebraically closed of characteristic not equal to two. Let $f\in(x_{1},\cdots,x_{n})^{2}$, and set $R=S/(f)$ and $R^{\sharp\sharp}=S[[x_{1},\cdots,x_{n+1},z_{1},z_{2}]]/(f+z_{1}^{2}+z_{2}^{2})$. There is a triangulated equivalence
\[
\stgf{R}\simeq \stgf{R^{\sharp\sharp}}
\] 
which, when restricted to compact objects, agrees with Kn\"{o}rrer's equivalence.  
\end{thm*} 
We note that to prove the above equivalence we use Kn\"{o}rrer's, so the result extends, rather than recovers, the classic result. As an application of this theorem, we show in \cref{hypziegler} there is an extremely close relationship between Ziegler spectra of $\msf{GFlat}(R)$ and $\msf{GFlat}(R^{\sharp\sharp})$. 

We then turn our attention to examining adjoint functors to the familiar change of rings functor  \[
R\otimes_{R^{\sharp\sharp}}-\colon\stgf{R^{\sharp\sharp}}\to\stgf{R},\] which, we note, is by no means an equivalence even though the categories are themselves equivalent. It is shown in \cref{ambidextrous} this functor admits an ambidextrous adjunction: it has a left and right adjoint which coincide. 

Finally, in \cref{sec:flatover} we focus on the specific case when $R\to S$ is a ring homomorphism such that $S$ is finitely generated projective over $R$. In this setting $S\otimes_{R}-$ admits both a left and right adjoint. We investigate when these adjoints preserve Gorenstein flat cotorsion modules and induce adjoints on the stable categories. Following this, we show that for ring epimorphisms the endofunctor $S\otimes_{R}-$ yields a recollement on the stable categories. More specifically,

\begin{thm*}[{\ref{epirecollement}}]
Let $R\to S$ be an epimorphism of right coherent rings such that $S$ is finitely presented projective over $R$. Then, $\msf{Im}$, the image of $S\otimes_{R}-\colon\gfc{R}\to\gfc{S}$, is a Frobenius category. Moreover, $S\otimes_{R}-\colon\stgf{R}\to\stgf{S}$ is a Bousfield localisation on $\stgf{R}$. In particular, there is a recollement of triangulated categories
\[
\begin{tikzcd}[column sep= 1in]
\msf{C} \arrow[r, "\mrm{inc}" description ] \arrow[r, leftarrow, shift left = 2ex, "\lambda"] \arrow[r, leftarrow, shift left = -2ex, "\rho", swap] &  \stgf{R} \arrow[r, "S\otimes_{R}-" description ] \arrow[r, leftarrow, shift left = 2ex, "S^{*}\otimes_{S}-"] \arrow[r, leftarrow, shift left = -2ex, "\msf{res}", swap] & \ul{\msf{Im}}
\end{tikzcd}
\]
where $\msf{C}=\{X\in\stgf{R}:S\otimes_{R}X=0\}$.
\end{thm*}

\section{Preliminaries}\label{recollections}

In this section the relevant background material for most of the document shall be recalled. This will encompass three main areas: purity, Gorenstein modules, and some properties of triangulated categories. However, before this, some notation is fixed.

Let $R$ be a ring. An $R$-module will always mean a left $R$-module, and the category of $R$-modules shall be denoted $\Mod{R}$. Right $R$-modules shall be referred to explicitly as such, and in general will be identified with modules over $R^{\op}$; hence $\Mod{R^{\op}}$ denotes the category of right $R$-modules. The full subcategory of $R$-modules which admit projective resolutions of finitely generated projective modules shall be denoted $\mod{R}$. For brevity, the category of abelian groups is denoted $\ab$.

\subsection{Purity in finitely accessible categories}

\begin{chunk}
Let $\A$ be an additive category with direct limits. An object $A\in\A$ is \ti{finitely presented} if the functor $\Hom_{\A}(A,-)\colon \A\to\ab$ preserves direct limits. Recall that $\A$ is \ti{finitely accessible} if the category $\msf{fp}(\A)$ of finitely presented objects in $\A$ is skeletally small and $\A=\rlim\msf{fp}(\A)$.
\end{chunk}

\begin{chunk}\label{pureexactsequence}
In any finitely accessible category $\A$ there is a notion of purity. As illustrated in \cite[Theorem 5.2]{dac}, a short exact sequence $0\to X\to Y\to Z\to 0$ in $\A$ is \ti{pure exact} provided any of the following equivalent conditions hold: 
\begin{enumerate}
\item the induced sequence
$0\to\Hom_{\A}(A,X)\to\Hom_{\A}(A,Y)\to\Hom_{\A}(A,Z)\to 0$ 
is exact in $\ab$ for every $A\in\msf{fp}(\A)$;
\item there is a directed system $\mc{S}=(0\to A_{i}\to B_{i}\to C_{i}\to 0)_{I}$ of split exact sequences in $\A$ such that $\rlim\mc{S}=(0\to A\to B\to C\to 0)$.
\end{enumerate} 
If $\A=\Mod{R}$ for a ring $R$, then the above conditions are further equivalent to
\begin{enumerate}[resume]
\item for any $M\in\Mod{R^{\op}}$, the induced sequence $0\to M\otimes_{R}X\to M\otimes_{R}Y\to M\otimes_{R}Z\to 0$ is exact.
\end{enumerate}

Given a pure exact sequence $0\to X\to Y\to Z\to 0$, we say that $X$ is a \ti{pure subobject} of $Y$, and that $Z$ is a \ti{pure quotient} of $Y$. The maps $X\to Y$ and $Y\to Z$ are a \ti{pure monomorphism} and \ti{pure epimorphism}, respectively. When $\A=\Mod{R}$, we shall just say that $X$ is a pure submodule of $Y$.
\end{chunk}

\begin{chunk}\label{pureinjectivemodule}
Let $\A$ be a finitely accessible category, and assume that $\A$ also has products. The objects which are injective with respect to the pure exact sequences shall take a central role throughout this document. Recall that if $A\in\A$, then for any set $I$ we let $A^{I}=\prod_{I}A$.

An object $X\in\A$ is \ti{pure injective} provided any of the following equivalent conditions hold:
\begin{enumerate}
\item $0\to\Hom_{R}(N,X)\to\Hom_{R}(M,X)\to\Hom_{R}(L,X)\to 0$ is exact in $\ab$ for any pure exact sequence $0\to L\to M\to N\to 0$;
\item for any set $I$, there is a map $f:X^{I}\to X$ such that $f\circ\lambda_{i}=\msf{Id}_{X}$, where $\lambda_{i}:X\to X^{I}$ is the canonical embedding;
\item any pure monomorphism $X\to Y$ is split.
\end{enumerate}
The equivalence of the first and third items can be found at \cite[Theorem 5.4]{dac}, while the equivalence with the second is discussed prior to \cite[Definition 5.1]{saorstov}.

The class of pure injective objects in $\A$ shall be denoted $\msf{Pinj}(\A)$. While this is a proper class, there is only a set of isomorphism classes of indecomposable pure injective objects of $\A$, and this set shall be denoted $\msf{pinj}(\A)$. The fact there is such a set can be found at \cite[Corollary 4.3.38]{psl}. This set can be equipped with a topology which we now provide the background for.
\end{chunk}

\begin{chunk}
For a finitely accessible category with products $\A$, we let $\Mod{\msf{fp}(\A)}$ denote the category of left $\msf{fp}(\A)$ modules, which is just the category of additive functors $\msf{fp}(\A)\to\ab$. This is, by \cite[Theorem 6.1]{dac}, a locally coherent Grothendieck category, i.e. a finitely accessible Grothendieck category whose finitely presented objects form an abelian category; we let $\mod{\msf{fp}(\A)}$ denote this abelian category.

An object $f\in\Mod{\msf{fp}(\A)}$ is finitely presented if and only if it has a presentation of the form
\[
\Hom_{\A}(A,-)\vert_{\msf{fp}(\A)}\to\Hom_{\A}(B,-)\vert_{\msf{fp}(\A)}\to f\to 0
\]
where $A,B\in\msf{fp}(\A)$. Such a finitely presented object $f$ extends uniquely to a functor $\bar{f}\colon\A\to\ab$ which preserves direct limits and products: it is the functor with the presentation
\[
\Hom_{\A}(A,-)\to\Hom_{\A}(B,-)\to\bar{f}\to 0
\]
in $(\A,\ab)$.
\end{chunk}

\begin{chunk}\label{definablesubcategories}
A subcategory $\mc{D}\subset\A$, where $\A$ is a finitely accessible category with products, is said to be \ti{definable} if there is a set $\mf{X}\subseteq\mod{\msf{fp}(\A)}$ such that
\[
\mc{D}=\{A\in\A:\bar{f}A=0 \t{ for all }f\in\mf{X}\}.
\]
It transpires that $\mc{D}$ being definable is equivalent to it being closed, in $\A$, under pure subobjects, products and direct limits, see \cite[Theorem 3.4.7]{psl}. 

In fact, the relationship between definable subcategories of $\A$ and $\mod{\msf{fp}(\A)}$ runs deeper. If $\mc{D}\subset\A$ is definable, then the set 
\[
\mathscr{S}(\mc{D})=\{f\in\mod{\msf{fp}(\A)}:\bar{f}X = 0 \t{ for all }X\in\mc{D}\},
\]
is a Serre subcategory of $\mod{\msf{fp}(\A)}$. Conversely, if $\mc{S}$ is a Serre subcategory of $\mod{\msf{fp}(\A)}$, then 
\[
\mathscr{D}(\mc{S})=\{X\in\A:\bar{f}X=0 \t{ for all }f\in\mc{S}\}
\]
is a definable subcategory of $\A$. The assignments $\mc{D}\mapsto\mathscr{S}(\mc{D})$ and $\mc{S}\mapsto\mathscr{D}(\mc{S})$ provide mutually inverse order reversing bijections between definable subcategories of $\A$ and Serre subcategories of $\mod{\msf{fp}(\A)}$. See \cite{kredc} and \cite{dac} for more details.
\end{chunk}

\begin{ex}\label{flatisdefinable}
Let us give an example of a definable subcategory that shall appear throughout. Recall that a ring $R$ is \ti{right coherent} if the category $\mod{R}$ is abelian, which is equivalent to every finitely generated right ideal of $R$ being finitely presented. We consider the class $\msf{Flat}(R)$ of flat $R$-modules. 

By definition, a module $M$ is in $\msf{Flat}(R)$ if and only if $\t{Tor}_{R}^{1}(R/I,M)=0$ for all right ideals  $R\supseteq I$. In particular, we see that $\msf{Flat}(R)$ is definable if and only if $\t{Tor}_{R}^{1}(R/I,-)\in\msf{fp}(\mod{R},\ab)$ is finitely presented. It is shown in \cite[Theorem 10.2.36]{psl} that this occurs whenever each $R/I$ has a presentation of the form $P_{2}\to P_{1}\to P_{0}\to R/I\to 0$, where $P_{i}$ is finitely generated. Yet this is clearly equivalent to $I$ being finitely presented for every right ideal $I$. 

In other words, $\msf{Flat}(R)$ is definable if and only if $R$ is right coherent. Of course, there is a more direct proof of this fact: $\msf{Flat}(R)$ is always closed under pure subobjects and direct limits, so one only needs to know when it is closed products, and this is known to happen whenever $R$ is right coherent by \cite[\S2]{chase}.
\end{ex}

\begin{chunk}\label{zieglerspectrumdef}
Given a finitely accessible category with products, $\A$, we may now define the topology on $\msf{pinj}(\A)$. Say a subset $\msf{X}\subseteq\msf{pinj}(\A)$ is \ti{Ziegler closed} if there is a definable subcategory $\mc{D}\subseteq\A$ such that $\mc{D}\cap\msf{pinj}(\A)=\msf{X}$. One need not worry that different definable subcategories yield the same closed set, for if $\mc{D}_{1},\mc{D}_{2}\subseteq\A$ are definable subcategories, then $\mc{D}_{1}=\mc{D}_{2}$ if and only if $\mc{D}_{1}\cap\msf{pinj}(\A)=\mc{D}_{2}\cap\msf{pinj}(\A)$, see \cite[Corollary 5.1.5]{psl}. In particular, definable subcategories are uniquely determined by the indecomposable pure injective objects they contain.

The \ti{Ziegler spectrum} is defined to be the topological space whose underlying set is $\msf{pinj}(\A)$, and whose closed sets are given by the Ziegler closed subsets. This space is denoted by $\msf{Zg}(\A)$. If $\mc{D}\subset\A$ is a definable subcategory, we let $\msf{Zg}(\mc{D})$ be the topological space $\msf{pinj}(\A)\cap\mc{D}$ equipped with the subspace topology.
\end{chunk}

\begin{chunk}\label{definablefunctor}
Having considered definable subcategories in isolation, we now introduce the right notion of a functor between them. Suppose that $\A$ and $\B$ are finitely accessible categories with products, and $\mc{C}\subset\A$ and $\mc{D}\subset\B$ are definable subcategories.

An additive functor $F\colon \mc{C}\to\mc{D}$ is \ti{definable} if it preserves direct limits and direct products. By the definition of a pure exact sequence in \cref{pureexactsequence}, we see that $F$ preserving direct limits means that it preserves pure exact sequences. Moreover, from the definition in \cref{pureinjectivemodule}, we see that $F$ preserves pure injective objects.

It is also not difficult to see that, if $\mc{C}'\subseteq \mc{C}$ and $\mc{D}'\subseteq\mc{D}$ are definable subcategories, then $F^{-1}(\mc{D}')=\{X\in\mc{C}:FX\in\mc{D}'\}$ is definable, as is $\msf{pure}(F\mc{C}')\subseteq\mc{D}$, the closure of the image of $F$ under pure submodules, see \cite[p.54]{dac} for both.
\end{chunk}

\begin{ex}\label{extensionflatexample}
As the latter part of this work focuses on ascent properties in relation to purity, we shall predominantly consider the functor $S\otimes_{R}-\colon\Mod{R}\to\Mod{S}$, where $S$ is a ring that is also a right $R$-module. 

Given a collection $\{M_{i}\}_{I}$ of $R$-modules, there is a canonical map
\[
\phi\colon S\otimes_{R}\prod_{I}M_{i}\to\prod_{i}S\otimes_{R}M_{i},
\]
which is an isomorphism if and only if $S$ is finitely presented over $R$. 

Since $S\otimes_{R}-$ always preserves direct limits and sends projective $R$-modules to projective $S$-modules, one may wonder if its restriction to $\Flat{R}\to\Flat{S}$ is definable whenever $R$ and $S$ are right coherent rings with a weaker assumption on $S$ than being in $\mod{R^{\op}}$. This is not possible. For, if $\phi$ is an isomorphism when each $M_{i}$ is flat, it must be an epimorphism when $M_{i}=R$ for each $i\in I$. But \cite[Lemma 3.8]{gt} states that this can only happen when $S$ is finitely generated over $R$. In this case, by \cite[Theorem 1]{goodearl}, the identity on $S$ factors through a finitely presented right $R$-module, and thus, as $R$ was assumed to be right coherent, $S$ is itself finitely presented over $R$.
\end{ex}

\begin{chunk}\label{maponziegler}
Given a definable functor $F\colon\mc{C}\to\mc{D}$ between two definable subcategories (of two finitely accessible categories with products), we need not immediately obtain an induced map of topological spaces $\msf{Zg}(\mc{C})\to\msf{Zg}(\mc{D})$, since $F$ need not preserve \ti{indecomposable} pure injective modules. The condition which ensures that indecomposability is preserved is that $F$ is full on pure injectives, that is for any pure injective objects $X,Y\in\msf{Pinj}(\mc{C})$, the natural map $\Hom_{\mc{C}}(X,Y)\to\Hom_{\mc{D}}(FX,FY)$ is surjective. 

With that condition satisfied, not only does one obtain a map of sets $\msf{Zg}(\mc{C})\to\msf{Zg}(\mc{D})$, but this is a closed and continuous map of topological spaces, which induces a homeomorphism $\msf{Zg}(\mc{C})\setminus \msf{K}\to \msf{Zg}(\mc{D})\cap\msf{pure}(F\mc{C})$, where $\msf{K}=\{X\in\msf{Zg}(\mc{C}):FX=0\}$. Proofs of these statements can be found in \cite[Chapter 15]{dac}.
\end{chunk}

\subsection{Approximations, cotorsion pairs, and Gorenstein modules}

\begin{chunk}
We now recall the notions of (pre)covers and (pre)envelopes. Suppose $\mc{F}\subseteq\Mod{R}$ is a class of modules, and $M\in\Mod{R}$ is an arbitrary $R$-module. A morphism $f\colon F\to M$, where $F\in\mc{F}$ is said to be an $\mc{F}$-\ti{precover} if for any other map $f'\colon F'\to M$, with $F'\in\mc{F}$, there is a map $\alpha\colon F'\to F$ such that $f'=f\circ \alpha$. An $\mc{F}$-precover $f\colon F\to M$ is a \ti{cover} if, whenever $g\in\t{End}_{R}(F)$ is such that $f\circ g=f$, then $g$ is an automorphism. The class $\mc{F}$ is called (pre)covering provided every module admits an $\mc{F}$-(pre)cover. The dual notions are called pre-envelopes and envelopes.
\end{chunk}

\begin{ex}\label{definableclosedunderpienvelopes}
Some of the classes we have already introduced admit approximations. Let $\A$ be a finitely accessible category with products and suppose $\mc{D}\subseteq\A$ is a definable subcategory. Then $\mc{D}$ is pre-enveloping and covering in $\A$. A proof of the former can be found at \cite[Proposition 3.4.42]{psl}, while the latter was proved at \cite[Theorem 2.6]{cpt}. 

The class $\msf{Pinj}(\A)$ is enveloping. This was illustrated originally for pure injective modules over a ring at \cite[Proposition 6]{warfield}, and a proof of the more general statement can be found at \cite[Proposition 12.2.8]{krbook}. 

Given an object $M$, we let $PE(M)$ denote its pure injective envelope. A fundamental property is that the enveloping map $M\to PE(M)$ is, in fact, a pure monomorphism; in other words, every module is a pure subobject of its pure injective envelope. Pure injective envelopes behave well with respect to definable subcategories: $M\in\mc{D}$ if and only if $PE(M)\in\mc{D}$, as shown at \cite[Theorem 4.3.21]{psl}.
\end{ex}

\begin{chunk}\label{cotorsionpair}
Another frequent source of approximations are cotorsion pairs. Since we shall only be concerned with cotorsion pairs in $\Mod{R}$, we give the definition in that setting, although it is worth noting the concept exists much more generally.

Given a class $\mc{A}\subset\Mod{R}$, we define
\[
\mc{A}^{\perp}=\{X\in\Mod{R}:\t{Ext}_{R}^{1}(A,X)=0 \t{ for all }A\in\mc{A}\}
\]
and
\[
^{\perp}\mc{A}=\{X\in\Mod{R}:\t{Ext}_{R}^{1}(X,A)=0 \t{ for all }A\in\mc{A}\}
\]
to be the right and left Ext-orthogonal classes to $\mc{A}$, respectively. A pair of classes $(\mc{A},\mc{B})$ of $R$-modules is called a \ti{cotorsion pair} if $\mc{A}^{\perp}=\mc{B}$ and $^{\perp}\mc{B}=\mc{A}$. 

We say that a class $\mc{A}$ is \ti{self-orthogonal} if $\mc{A}=\mc{A}\cap\mc{A}^{\perp}=\,^{\perp}\mc{A}\cap\mc{A}$. Note that if $(\mc{A},\mc{B})$ is a cotorsion pair, then $\mc{A}\cap\mc{B}$ is a self-orthogonal class.

As mentioned, cotorsion pairs frequently provide well behaved approximations. Given a class $\mc{F}$, an $\mc{F}$-precover $\alpha\colon F\to M$ is \ti{special} if $\alpha$ is surjective and $\msf{Ker}(f)\in\mc{F}^{\perp}$, and the dual notion defines a \ti{special pre-envelope}. For a cotorsion pair $(\mc{A},\mc{B})$, the existence of approximations by $\mc{A}$ and $\mc{B}$ are closely linked through \ti{Salce's lemma} (see \cite[Lemma 5.20]{gt}), which states the following are equivalent:
\begin{enumerate}
\item $\mc{A}$ is special pre-covering,
\item $\mc{B}$ is special pre-enveloping.
\end{enumerate}
In this case, the cotorsion pair $(\mc{A},\mc{B})$ is called \ti{complete}. When the classes of the cotorsion pair $(\mc{A},\mc{B})$ are covering and enveloping, respectively, then the cotorsion pair is called \ti{perfect}; and any perfect cotorsion pair is complete.
\end{chunk}

\begin{ex}\label{Flatcotorsionpair}
Over any ring $R$, there is a cotorsion pair $(\Flat{R},\msf{Cot}(R))$, where the modules in $\msf{Cot}(R)$ are called \ti{cotorsion} modules. In \cite{fcc}, it was shown that $(\msf{Flat}(R),\msf{Cot}(R))$ is a perfect and hereditary cotorsion pair. In particular, special flat covers always exist, as do special cotorsion envelopes. As a consequence of this, over any ring the cotorsion envelope of a flat module is flat, while the flat cover of any cotorsion module is cotorsion. Indeed, if $F$ is flat and $0\to F\to C(F)\to Q\to 0$ is a special cotorsion envelope of $F$, then as $Q$ is flat and the flat modules are extension closed, we have that $C(F)$ is also flat; the dual statement is proved in the same way. The intersection $\Flat{R}\cap\msf{Cot}(R)$ shall be of significant use throughout, and we denote it by $\msf{FlatCot}(R)$.

An observation that will be crucial to us is that any pure injective module is cotorsion. To see this, note that from the third condition of \cref{pureexactsequence}, any exact sequence ending in a flat module is pure exact, and every flat module appears as the end of a short exact sequence by considering a presentation. 

Consequently, by definition of being pure injective as in \cref{pureinjectivemodule}, we have that $\t{Ext}_{R}^{1}(F,X)=0$ for all $F\in\Flat{R}$ and $X\in\msf{Pinj}(R)$. Thus $\Flat{R}\cap\msf{Pinj}(R)\subseteq\msf{FlatCot}(R)$.

One may wonder when the reverse inclusion holds, that is when $\msf{FlatCot}(R)=\Flat{R}\cap\msf{Pinj}(R)$. This was shown to hold for right coherent rings in \cite[Lemma 3.2.3]{xu}. In more generality, it was shown in \cite[Theorem 3.2]{rothmaler} that $\msf{FlatCot}(R)=\msf{Flat}(R)\cap\msf{Pinj}(R)$ is equivalent to the property that $PE(F)\in\Flat{R}$ for all $F\in\Flat{R}$. Of course, this happens whenever $R$ is right coherent, by \cref{flatisdefinable} and \cref{definableclosedunderpienvelopes}, as in this setting $\Flat{R}$ is a definable subcategory of $\Mod{R}$.
\end{ex}

\begin{chunk}\label{totallyacyclicdef}
Let $\mc{A}$ be a self-orthogonal class of modules. In \cite{CET}, Gorenstein objects with respect to $\mc{A}$ were introduced, which recovers some classical constructions in Gorenstein homological algebra. To introduce these Gorenstein objects, the appropriate relative notion of a totally acyclic complex was given.

Let $\mc{A}\subset\Mod{R}$ be a self orthogonal class. A chain complex $T$ over $R$ is $\mc{A}$-\ti{totally acyclic} provided the following conditions hold:
\begin{enumerate}
\item $T$ is acyclic;
\item $T_{i}\in\mc{A}$ for every $i\in\Z$;
\item $\Hom_{R}(T,A)$ is acyclic for all $A\in\mc{A}$;
\item $\Hom_{R}(A,T)$ is acyclic for all $A\in\mc{A}$.
\end{enumerate}

This is not the definition as given in \cite[Definition 1.1]{CET}, rather it is a combination of \cite[Propositions 1.3 and 1.5]{CET}, which is sufficient to our goals, as we shall only consider self orthogonal classes.
\end{chunk}

\begin{ex}
There are two immediate self-orthogonal classes: $\msf{Proj}(R)$ of projective $R$-modules, and $\msf{Inj}(R)$ of injective $R$-modules. The class of $\msf{Proj}$-totally acyclic complexes has long been an object of study, and is usually referred to as a totally acyclic complex of projective modules. The same holds true for $\msf{Inj}(R)$-totally acyclic complexes. See, for example, \cite{ik}, \cite{krstable}, \cite{ops}, among many more. 

The class of $\msf{FlatCot}$-totally acyclic complexes shall play a crucial role throughout, but we shall wait before introducing it.
\end{ex}

\begin{chunk}\label{Agorenstein}
Let $\mc{A}\subset\Mod{R}$ be a self-orthogonal class. An $R$-module $M$ is $\mc{A}$-Gorenstein if $M=Z_{0}(T)$ for some $\mc{A}$-totally acyclic complex $T$.

The $\msf{Inj}(R)$-Gorenstein objects are nothing other than the usual Gorenstein injective $R$-modules, while the $\msf{Proj}(R)$-Gorenstein objects are the Gorenstein projective $R$-modules. We shall denote these classes by $\msf{GInj}(R)$ and $\msf{GProj}(R)$. These classes are frequently studied; information on them can be found in \cite{chris,dcmca,rha,iacob}, to name a few.
\end{chunk}

\begin{chunk}\label{gorensteinflat}
There is a third type of Gorenstein modules which cannot arise as a class of $\mc{A}$-totally acyclic complexes, in the sense of \cref{Agorenstein}, namely the Gorenstein flat modules.

Recall than an $R$-module $M$ is \ti{Gorenstein flat} if $M=Z_{0}(F)$ for an acyclic complex $F$ of flat $R$-modules such that $E\otimes_{R}F$ is acyclic for all $E\in\msf{Inj}(R^{\op})$. The full subcategory of Gorenstein flat modules is denoted $\msf{GFlat}(R)$.

An acyclic complex $F$ of flat $R$-modules such that $E\otimes_{R}F$ is acyclic for all $E\in\msf{Inj}(R^{\op})$ is frequently called an $\msf{F}$-\ti{totally acyclic complex}.
\end{chunk}

\begin{chunk}\label{GFlatcotorsionpair}
Let us make some observations on the class $\msf{GFlat}(R)$. Since, for any $F\in\Flat{R}$, the complex $0\to F\xrightarrow{\t{Id}}F\to 0$ is $\msf{F}$-totally acyclic, any flat module is Gorenstein flat. In many ways, $\msf{GFlat}(R)$ is a generalisation of $\Flat{R}$: for example, over any ring there is a complete hereditary cotorsion pair $(\msf{GFlat}(R),\msf{GCot}(R))$ by \cite[Corollary 4.12]{sarstov}, generalising \cite[Theorem 2.11]{ejlr}, which proved the existence of such a cotorsion pair whenever $R$ was right coherent. The modules in $\msf{GCot}(R)=\msf{GFlat}(R)^{\perp}$ are called \ti{Gorenstein cotorsion}. As illustrated at \cite[Corollary 4.12]{sarstov}, there is an equality $\msf{FlatCot}(R)=\msf{GFlat}(R)\cap\msf{GCot}(R)$. 
\end{chunk}

\begin{chunk}\label{gflatproducts}
Similar to the flat modules, the class $\msf{GFlat}(R)$ need not, in general, be closed under products; it remains an open question when this is the case. However, it was shown at \cite[Proposition 4.13]{sarstov} that $\msf{GFlat}(R)$ is closed under products if and only if it is a definable subcategory of $\Mod{R}$, in which case $R$ is necessarily right coherent. Since the focus on this paper is when the class $\msf{GFlat}(R)$ is definable, we may therefore always assume that $R$ is right coherent, in which case, of course, $\msf{Flat}(R)$ is also definable. In this setting, it follows from \cref{Flatcotorsionpair} and \cref{GFlatcotorsionpair} that any module that is simultaneously Gorenstein flat and Gorenstein cotorsion is necessarily pure injective. 
\end{chunk}

\begin{chunk}\label{gpgfrelations}
However, there is one striking difference between the flat and Gorenstein flat modules. Over any ring it is a classic result of Govorov and Lazard that $\Flat{R}$ is a finitely accessible category, with $\msf{fp}(\Flat{R})=\msf{proj}(R)$, the category of finitely generated projective $R$-modules. The Gorenstein analogue is not true: as illustrated in \cite{glbk,glhj}, there are rings such that that $\msf{GFlat}(R)\not=\rlim\msf{Gproj}(R)$. However, over Gorenstein rings it is the case that $\msf{GFlat}(R)$ is finitely accessible with $\msf{fp}(\msf{GFlat}(R))=\msf{Gproj}(R)$, as illustrated in \cite[Theorem 10.3.8]{rha}.

Furthermore, it is not currently known whether there is even always an inclusion $\msf{Gproj}(R)\subset\msf{GFlat}(R)$. Over many rings this inclusion is known, for example it holds for any right coherent ring such that every flat module has finite projective dimension, by \cite[Theorem 4]{iacob2}. In particular, this holds over any Gorenstein ring, or any Noetherian ring admitting a dualising complex. This is further discussed in \cref{gprojmodel}.
\end{chunk}

\begin{chunk}\label{gorflatcot}
If $F$ is an $\msf{F}$-totally acyclic complex, such that $F_{i}\in\msf{FlatCot}(R)$, then it was shown in \cite[Theorem 1.3]{bce} that $Z_{0}(M)$ is a cotorsion module, as well as a Gorenstein flat module. In fact, \cite[Theorem 4.4]{CET} shows that, whenever $R$ is a right coherent ring, a complex $F$ is an $\msf{F}$-totally acyclic complex that is degreewise cotorsion if and only if it is a $\msf{FlatCot}(R)$-totally acyclic complex, in the sense of \cref{totallyacyclicdef}.

Moreover, it was shown in \cite[Theorem 5.2]{CET} that, with the assumption of $R$ being right coherent, an $R$-module is Gorenstein flat and cotorsion if and only if it is $\msf{FlatCot}(R)$-Gorenstein. We shall denote this class by $\gfc{R}$, and call a $\msf{FlatCot}(R)$-totally acyclic complex a \ti{totally acyclic complex of flat cotorsion modules}.

Before giving more details about $\gfc{R}$, we shall introduce some of the tools from triangulated categories.
\end{chunk}

\subsection{Frobenius categories, triangulated categories, and their purity}

\begin{chunk}\label{frobeniusdef}
Recall that any extension closed subcategory $\msf{E}$ of $R$-modules has the structure of an exact category, in the sense of Quillen, see \cite{buhler} for details. The short exact sequences in $\msf{E}$ are just the short exact sequences of $R$-modules that are termwise in $\msf{E}$. If $0\to A\to B\to C\to 0$ is a short exact sequence in $\msf{E}$, we call the map $A\to B$ an inflation and the map $B\to C$ a deflation.

An object $E\in\msf{E}$ is \ti{injective} if any inflation $E\to X$ splits, while an object $P\in \msf{E}$ is \ti{projective} if any deflation $Y\to P$ splits. The category $\msf{E}$ is said to have \ti{enough injectives} (resp. \ti{enough projectives}) if for any $X\in\msf{E}$ there is an inflation $X\to E$ with $E$ injective, respectively a deflation $P\to X$ with $P$ projective.

An exact category $\msf{E}$ is \ti{Frobenius} if it has enough projective and injective objects, and these two classes coincide. 

Given a Frobenius category $\msf{E}$, one can form its stable category $\ul{\msf{E}}$, which is a triangulated category. This is the category whose objects are the same as those in $\msf{E}$, but
\[
\underline{\Hom}_{\msf{E}}(X,Y):=\Hom_{\ul{\msf{E}}}(X,Y)=\Hom_{\msf{E}}(X,Y)/\{f\colon X\to Y:\t{ $f$ factors through a projective object}\}.
\]
The triangulated structure on $\msf{E}$ is straightforward to describe: the shift $\Sigma X$ of $X\in\ul{\msf{E}}$ is given by the first cosyzygy, while $X\to Y\to Z\to\Sigma X$ is a triangle in $\ul{\msf{E}}$ if and only if $0\to X\to Y\to Z\to 0$ is a short exact sequence in $\msf{E}$. These two facts are proved, and described in more detail, in \cite[\S 3.3]{krbook}.

\end{chunk}

\begin{chunk}\label{gfcfrobenius}
A specific case of \cite[Theorem 2.11]{CET} shows that, whenever $\mc{A}\subset\Mod{R}$ is a self-orthogonal class of modules, then the class of $\mc{A}$-Gorenstein objects is a Frobenius category, whose projective-injective objects are precisely the objects of $\mc{A}$.

This recovers some well known results: it was shown in \cite[\S 7]{krstable} that $\msf{GInj}(R)$ was a Frobenius category, while it was shown by Buchweitz in \cite{buchweitz} that $\msf{Gproj}(R)$ was a Frobenius category; it has also long been known that $\msf{GProj}(R)$ is a Frobenius category.

However, $\msf{GFlat}(R)$ is seldom a Frobenius category. This, as shown in \cite[Theorem 4.5]{CET}, is equivalent to $R$ being left perfect. However, the same theorem shows that, over any ring, the class $\gfc{R}$ is Frobenius, with injective objects being $\msf{FlatCot}(R)$. 

In particular, when $R$ is right coherent, we have, as discussed in \cref{gorflatcot}, that $\gfc{R}$ is the class of $\msf{FlatCot}$-Gorenstein objects. Since we shall be interested in rings where $\msf{GFlat}(R)$ is definable, which, by \cref{gflatproducts} necessitates $R$ being right coherent, we may always assume that $\gfc{R}$ is this class of $\msf{FlatCot}$-Gorenstein objects.
\end{chunk}

\begin{chunk}
The stable categories of the Frobenius categories of $\mc{A}$-Gorenstein objects can be explicitly given. In \cite[Theorem 3.8]{CET} it is shown that there are triangulated equivalences
\[
\begin{tikzcd}[column sep=2cm]
\ul{\mc{A}\t{-}\msf{Gor}} \arrow[r, shift left =1ex, "\mrm{T}"] \arrow[r, shift right=1ex, leftarrow, swap, "Z_{0}"] & \msf{K}_{\t{tac}}(\mc{A})
\end{tikzcd}
\]
between the stable category of $\mc{A}$-Gorenstein objects and the homotopy category of $\mc{A}$-totally acyclic complexes. Here $T$ is the functor sending an object to the totally acyclic complex it is a cycle in.

As a specific case, the above gives a triangulated equivalence $\stgf{R}\to\msf{K}_{\t{tac}}(\msf{FlatCot}(R))$, which shall be of great utility throughout.
\end{chunk}

\begin{chunk}\label{purecaveat}
We shall be interested in purity in these stable categories. However, there is a caveat here, namely that purity in triangulated categories is only understood when the categories are standard well-generated (see \cite{saorstov}), with substantially more understood when they are compactly generated. Since homotopy categories are frequently very far from being well generated, see \cite{stov} for some examples, it will not be possible to work in complete generality.

The main obstacle is in defining the concept of a pure triangle, as the notion of a pure injective object exists in any category with products - just use the second equivalent condition of \cref{pureinjectivemodule}.

To remedy this, we shall, when discussing purity in homotopy categories, impose the restriction that they are compactly generated, the definition of which we now recall.
\end{chunk}

\begin{chunk}\label{wellgenerated}
Let $\T$ be a triangulated category with coproducts and $\kappa$ a regular cardinal. $\T$ is said to be $\kappa$-\emph{well generated} if there is a set of objects $\scr{S}$ such that
\begin{enumerate}
\item if $\Hom_{\T}(S,X)=0$ for all $S\in\scr{S}$, then $X=0$;
\item every $S\in\scr{S}$ is $\kappa$-small, that is if $f\colon S\to \oplus_{I}X_{i}$ is any morphism, then there is a subset $J\subseteq I$ with $\vert J\vert <\kappa$ such that $f$ factors as $S\to\oplus_{j\in J}X_{j}\to\oplus_{I}X_{i}$. Here the last map is the canonical split embedding;
\item for any morphism $f\colon S\to \oplus_{i}X_{i}$ with $S\in\scr{S}$, there are maps $f_{i}\colon S_{i}\to X_{i}$, with $S_{i}\in\scr{S}$ such that $f$ factors as $S\to \oplus_{i}S_{i}\to \oplus_{i}X_{i}$, where the last map is $\oplus_{I}f_{i}$.
\end{enumerate}
An $\aleph_{0}$-well generated triangulated category is called \emph{compactly generated}.  $\T$ is said to be \emph{well generated} if it is $\kappa$-well generated for some $\kappa$.
\end{chunk}

\begin{chunk}\label{puritytriangulateddef}
Let us give a brief description of purity in compactly generated triangulated categories; a more detailed exposition can be found in \cite{birdwilliamsonhomological}.

Let $\T$ be a compactly generated triangulated category with compact objects $\T^{\c}$. The main tool for understanding purity in $\T$ is the \ti{restricted Yoneda embedding}
\[
\y\colon\T\to\Mod{(\T^{\c})^{\op}}=((\T^{\c})^{\op},\ab)
\]
\[
X\mapsto \Hom_{\T}(-,X)\vert_{\T^{\c}}.
\]
Following \cite{krsmash}, a triangle $X\to Y\to Z\to \Sigma X$ in $\T$ is said to be \ti{pure} if the sequence $0\to\y X\to \y Y \to \y Z\to 0$ is exact in $\Mod{(\T^{\c})^{\op}}$. In this case the map $X\to Y$ is a \ti{pure monomorphism} and $Y\to Z$ is a pure epimorphism. Note that a triangle being pure is equivalent to the connecting map $Z\to \Sigma X$ being \ti{phantom}, that is $\Hom_{\T}(C,f)=0$ for all $C\in\T^{\c}$.

An object $X\in\T$ is is \ti{pure injective} if every pure monomorphism with domain $X$ splits; this is equivalent to $X$ satisfying the product criterion of \cref{pureinjectivemodule}. As shown in \cite[Theorem 1.8]{krsmash}, an object in $X$ is pure injective if and only if $\y X$ is an injective module in $\Mod{(\T^{\c})^{\op}}$; in particular there is again only a set of indecomposable pure injective objects in $\T$, which is denoted $\msf{pinj}(\T)$.
\end{chunk}

\begin{chunk}
Similarly to finitely accessible categories with products, we may also consider definable subcategories in compactly generated triangulated categories. A subcategory $\mc{D}\subset\T$ is \ti{definable} if there is a Serre subcategory $\mc{S}$ of $\mod{\T^{\c}}$ such that 
\[
\mc{D}=\{X\in\T:\bar{f}X=0\t{ for all }f\in\mc{S}\},
\]
where $\bar{f}$ is extended to a functor $\T\to \ab$ in the same way as in \cref{definablesubcategories}. Whenever $\T$ admits an enhancement, which is the case throughout this paper, a class $\mc{D}\subseteq\T$ is definable if and only if $\T$ is closed under pure subobjects, pure quotients and products, or equivalently filtered homotopy colimits, pure subobjects and products; this was proved in \cite[Theorem 3.11]{laking}.

Using definable subcategories of $\T$ and pure injective objects, we may also define the Ziegler spectrum of $\T$, which is denoted $\msf{Zg}(\T)$. The objects of $\msf{Zg}(\T)$ are the indecomposable pure injective objects of $\T$, while the closed sets are those of the form $\msf{pinj}(\T)\cap\mc{D}$, where $\mc{D}$ is a definable subcategory of $\T$.

In fact, purity in $\T$ can be completely described in terms of purity in the finitely accessible category $\Flat{(\T^{\c})^{\op}}$ of cohomological functors $(\T^{\c})^{\op}\to\ab$. This category can be realised as the ind-completion of $\T^{\c}$, viewed as an additive category. This was formalised in \cite[Lemma 2.9]{birdwilliamsonhomological}, but essentially $\y$ gives a bijection between pure triangles in $\T$ and pure exact sequences in $\Flat{(\T^{\c})^{\op}}$, as well as definable categories and pure injective objects. In particular, $\y$ induces a homeomorphism between $\msf{Zg}(\T)$ and $\msf{Zg}(\Flat{(\T^{\c})^{\op}})$.

\end{chunk}

\begin{chunk}\label{brownrep}
One final tool of triangulated categories that shall be used is Brown representability. This states that if a triangulated category $\T$ has arbitrary coproducts, then a functor $F\colon\T^{\op}\to \ab$ which is homological and sends coproducts to products is of the form $F\simeq \Hom_{\T}(-,X)$ for some $X\in\T$. This statement has the following consequences:
\begin{enumerate}
\item\label{brcoprod} $F\colon \T\to\U$ preserves coproducts if and only if it admits a right adjoint;
\item\label{brprod} $F\colon\T\to\U$ preserves products if and only if it admits a left adjoint.
\end{enumerate}
We note that Brown representability does not hold for all triangulated categories, but it does hold for compactly generated triangulated categories, see \cite[Proposition 5.3.1]{krloc}. 

We also note that (\ref{brcoprod}) holds for arbitrary well generated triangulated categories, by \cite[Theorem 8.4.4]{Neemantri}, while (\ref{brprod}) holds for homotopy categories of presentable stable $\infty$-categories, see the discussion at \cite[p. 17]{efimov}.
\end{chunk}

\section{Purity in the stable category of Gorenstein flat cotorsion modules}\label{purestructureofstable}

For this section, we shall assume that $\msf{GFlat}(R)$ is a definable subcategory of $\Mod{R}$. As mentioned in \cref{gflatproducts}, this forces $R$ to be a right coherent ring, hence $\Flat{R}$ is also definable. Moreover, we may identify $\gfc{R}$ with the $\msf{FlatCot}$-Gorenstein objects, hence $\gfc{R}$ is a Frobenius category with projective-injective objects $\msf{FlatCot}(R)$.

\subsection{Products and coproducts in the stable category}

\begin{chunk}
The first thing we shall investigate is the relationship between pure injectivity in $\msf{GFlat}(R)$ and $\stgf{R}$. By the discussion at \cref{Flatcotorsionpair}, any pure injective Gorenstein flat module is cotorsion, and therefore automatically lies in $\msf{GFlatCot}(R)$.

Thus understanding pure injective objects in $\msf{GFlat}(R)$ is equivalent to understanding pure injective objects in $\msf{GFlatCot}(R)$. We now show that this is equivalent to understanding pure injective objects in the stable category. First, note that cotorsion and pure injective objects are always closed under products, and with the assumptions imposed on $R$, we have that both $\gfc{R}$ and $\msf{FlatCot}(R)$ are product closed.  
\end{chunk}

\begin{lem}\label{piinst}
Suppose that $R$ is a right coherent ring such that $\msf{GFlat}(R)$ is closed under products. Then an object is pure injective in $\msf{GFlatCot}(R)$ if and only if it is pure injective in $\stgf{R}$.
\end{lem}

\begin{proof}
The localisation functor $\msf{GFlatCot}(R)\to\stgf{R}$ preserves products, hence, by the definition of pure injective module given at \cref{pureinjectivemodule}, the functor preserves pure injective objects. For the reverse implication, suppose that $X\in\stgf{R}$ is pure injective. Then, for any index set $I$, there is a morphism $f\colon X^{I}\to X$ such that $\msf{Id}_{X}=f\circ\lambda_{i}$. By definition of the morphisms in the stable category, there is a morphism $\tilde{f}\colon X^{I}\to X$ in $\gfc{R}$ such that $\msf{Id}_{X}-\tilde{f}\circ\lambda_{i}=\beta\circ\alpha$, where $\alpha\colon X\to F$ and $\beta\colon F\to X$ are maps with $F$ a flat cotorsion module. Since $F$ is an injective object in $\gfc{R}$, see \cref{gfcfrobenius}, and $\lambda_{i}$ is an inflation in $\gfc{R}$ for all $i$, there is a map $\gamma\colon X^{I}\to F$ such that $\alpha = \gamma\circ\lambda_{i}$. Then $\msf{Id}_{X}=(\tilde{f}+\beta\circ\gamma)\lambda_{i}$. Consequently, again by \cref{pureinjectivemodule}, we see that $X$ is pure injective in $\gfc{R}$.
\end{proof}

\begin{rem}
\cref{piinst} is a modification of \cite[Proposition 1.16]{krsmash}. In that setting, both products and coproducts are respected when passing to the stable module category.
\end{rem}

\begin{chunk}
As mentioned above, for \cref{piinst} we used the fact that products in $\gfc{R}$ agree with those in $\stgf{R}$, and they are the same as in $\Mod{R}$. One may wonder what happens with coproducts; here the situation is markedly different. Before stating the result, recall that $X\in\msf{Pinj}(R)$ is $\Sigma$-\ti{pure injective} if $X^{(I)}:=\oplus_{I}X$ is pure injective for any indexing set $I$; see \cite[\S 4.4.2]{psl}. 
\end{chunk}

\begin{lem}\label{coproductsperfect}
Let $R$ be a right coherent ring. Then the following are equivalent:
\begin{enumerate}
\item $\gfc{R}$ is closed under arbitrary coproducts in $\Mod{R}$;
\item $\msf{FlatCot}(R)$ is closed under arbitrary coproducts in $\Mod{R}$;
\item every flat $R$-module is cotorsion;
\item every $R$-module is cotorsion;
\item $R$ is a left perfect ring.
\end{enumerate}
\end{lem}
\begin{proof}
$(1)\Rightarrow (2)$ is straightforward.  For $(2)\Rightarrow (3)$, suppose $F\in\Flat{R}$, and consider $PE(F)\in\msf{FlatCot}(R)$ by \cref{Flatcotorsionpair}. By assumption, $PE(F)^{(I)}\in\msf{FlatCot}(R)$ for any set $I$, hence any coproduct of copies of $PE(F)$ is pure injective by \cref{Flatcotorsionpair}; in particular $PE(F)$ is $\Sigma$-pure injective. However, any pure submodule of a $\Sigma$-pure injective module is also $\Sigma$-pure injective by \cite[Proposition 4.4.12]{psl}, and thus $F$ is also $\Sigma$-pure injective and flat. In particular $F$ is pure injective, and thus is in $\msf{FlatCot}(R)$. 

For $(3)\Rightarrow(4)$, let $M\in\Mod{R}$ and consider the exact sequence $0\to C\to F(M)\to M\to 0$, where $F(M)\to M$ is the flat cover of $M$ and $C\in\msf{Cot}(R)$. This exists by \cref{Flatcotorsionpair}. By assumption, $F(M)$ is cotorsion; applying $\Hom_{R}(F,-)$ to this exact sequence for any $F\in\Flat{R}$, we obtain a monomorphism $\t{Ext}_{R}^{1}(F,M)\to\t{Ext}_{R}^{2}(F,C)$, and the latter is zero, as $C\in\msf{Cot}(\T)$. Thus $\t{Ext}_{R}^{1}(F,M)=0$ for any flat $R$-module $F$, hence $M\in\msf{Cot}(R)$.

$(4)\Rightarrow (5)$ is \cite[Proposition 3.2]{rothmaler}, once one makes the trivial observation that if every $R$-module is cotorsion, then $R$ is $\Sigma$-pure injective. Lastly $(5)\Rightarrow(1)$ is the contents of \cite[Theorem 4.5]{CET}, which states that under these assumptions, we have $\msf{GFlat}(R)=\gfc{R}$, and $\msf{GFlat}(R)$ is always closed under coproducts.
\end{proof}

\begin{chunk}\label{coproducts}
However, coproducts necessarily exist in $\stgf{R}$. This is because, over a right coherent ring, there is an abelian model structure on $\Mod{R}$ called the Gorenstein flat model structure, as introduced in \cite[Theorem 3.3]{gill}. This has the following properties:
\begin{enumerate}
\item the cofibrant objects are $\msf{GFlat}(R)$;
\item the fibrant objects are $\msf{Cot}(R)$;
\item the trivially cofibrant objects are $\msf{Flat}(R)$;
\item the trivially fibrant objects are $\msf{GCot}(R)$.
\end{enumerate}  
In particular, as shown in \cite[Corollary 3.4]{gill}, the homotopy category of the Gorenstein flat model structure is triangle equivalent to $\stgf{R}$, which then has arbitrary coproducts as all homotopy categories of model categories do.

More can be said. A result of Rosick\'{y}, see \cite[Theorem 12.21]{gillbook}, tells us that the stable category $\stgf{R}$ is a well generated triangulated category. Furthermore, both forms of Brown representability \cref{brownrep}.(\ref{brcoprod}) and \cref{brownrep}.(\ref{brprod}) hold in $\stgf{R}$ (see \cite[Proposition A.3.7.6]{htt} and the discussion in \cref{brownrep} for a justification). 

The coproducts in the Gorenstein flat model structure can be described as follows. The Gorenstein flat model structure is cofibrantly generated, by \cite[Corollary 4.12]{sarstov}, and thus we may assume taking fibrant replacements may be done functorially. The fibrant objects in the model structure are $\msf{Cot}(R)$, while the cofibrant objects are $\msf{GFlat}(R)$. Consequently, given $\{X_{i}\}_{I}\subseteq\gfc{R}$, the coproduct in $\gfc{R}$ is given by $\Phi(\oplus_{i} X_{i})$, where $\Phi$ is a fibrant replacement functor.

It is also possible to obtain the coproducts in $\stgf{R}$ through the homotopy category $\msf{K}(\Flat{R})$. It is clear that, whenever $R$ is a right coherent ring, the category $\msf{K}_{\t{tac}}(\msf{FlatCot}(R))$ is closed under products (as products in $\Mod{R}$ are exact), and Rosick\'{y}'s result gives that this category is well generated. As such, we may apply Brown representability \cref{brownrep} to the product preserving inclusion $\msf{K}_{\t{tac}}(\msf{FlatCot}(R))\to\msf{K}(\Flat{R})$, which thus admits a left adjoint. This left adjoint will create coproducts in $\msf{K}_{\t{tac}}(\msf{FlatCot}(R))$.
\end{chunk}

\begin{chunk}\label{gprojmodel}
If we assume that every flat $R$-module has finite projective dimension, we may also get coproducts in a different way, through the Gorenstein projective model structure, which is now described.

For an arbitrary ring $A$, it was shown in \cite[Theorem 8.5]{BGH} that there is a Gorenstein AC-projective model structure on $\Mod{A}$. For this structure, every object is fibrant and the cofibrant objects are the Gorenstein AC-projective modules. The pertinent fact for our purposes is that over coherent rings Gorenstein AC-projective modules coincide with the so called Ding projective modules, see \cite[\S 3.2]{gillding}. 

Now, in the case that $A$ is coherent such that every flat module has finite projective dimension, the Ding projective modules coincide with the usual Gorenstein projective modules. Indeed, \cite[Lemma 4]{iacob2} states that there is, for any ring $A$, an equality between Ding projective modules and Gorenstein projective modules if and only if $\msf{Flat}(A)\subseteq\msf{GProj}(A)^{\perp}$. Now, over any ring $A$, one has $M\in\msf{GProj}(A)^{\perp}$ whenever $M$ has finite projective or injective dimension by \cite[Lemma 2.1]{cfh}. 

We therefore see that if $A$ is a coherent ring such that every flat module has finite projective (or injective, in which case $A$ is just Gorenstein) dimension, then there is an equality between the Gorenstein projective modules and the Ding projective modules. With this holding, the homotopy category of this model structure is triangle equivalent to $\ul{\msf{GProj}}(R)$. Furthermore, we can deduce from \cite[Theorem 5.1]{estgil} that the identity functor provides a left Quillen equivalence between the Gorenstein projective and Gorenstein flat model structures, and thus there is a triangulated equivalence $\ul{\msf{GProj}}(R)\simeq \stgf{R}$.

As $\ul{\msf{GProj}}(R)$ is closed under coproducts inherited from $\Mod{R}$, we may pass these through the equivalence to endow $\stgf{R}$ with coproducts. These coproducts correspond to those induced by the triangulated equivalence between $\msf{K}_{\t{tac}}(\msf{FlatCot}(R))$ and $\msf{K}_{\t{tac}}(\msf{Proj}(R))$. The coproducts on the latter are just those of $\msf{K}(R)$. 
\end{chunk}

\subsection{Compactness, pure triangles and the Ziegler spectrum of the stable category}

While products are vital for pure injectivity in $\stgf{R}$, it is the coproducts which are fundamental to compactness, and compactness is indispensable to understand pure triangles and thus definability and the Ziegler spectrum. We now focus on these latter concepts.

\begin{chunk}
Having shown that an object is pure injective in $\msf{GFlat}(R)$ if and only if it is pure injective in $\stgf{R}$, one may wonder if the other component of the pure structure on $\msf{GFlat}(R)$, namely the pure exact sequences, can also be determined by the pure triangles in $\stgf{R}$, and vice versa.

As stated in \cref{purecaveat}, we must assume here that $\stgf{R}$ is compactly generated, yet this is not always the case. Part of the issue here is that, unless $R$ is left perfect, the category $\msf{K}(\Flat{R})$ is not even well generated, see \cite[Theorem 5.2]{stov}. However, over right coherent rings for which every flat left module has finite projective dimension, there is an equivalence $\msf{K}_{\t{tac}}(\msf{FlatCot}(R))\simeq\msf{K}_{\t{tac}}(\msf{Proj}(R))$, and the generation properties of the latter category are better understood.

For example \cite[Theorem 5.3]{ik} shows $\msf{K}_{\t{tac}}(\msf{Proj}(R))$ is compactly generated whenever $R$ is Noetherian and admits a dualising complex. As the existence of a dualising complex ensures that any flat module has finite projective dimension (see \cite[Proposition 3.4]{ik}), it follows that $\stgf{R}$ is also compactly generated in this setting. We may also appeal to the Gorenstein projective model structure for cases when $\ul{\msf{GProj}}(R)$ is compactly generated. By \cite[Theorem 8.42]{gillbook}, this occurs whenever the cotorsion pair $(\msf{GProj}(R),\msf{GProj}(R)^{\perp})$ is of finite type, that is $\msf{GProj}(R)^{\perp}=\mc{S}^{\perp}$ for some class of finitely presented $R$-modules $\mc{S}$.

With the assumption that $\stgf{R}$ is compactly generated, we shall let $\scr{C}$ denote the category of compact objects in $\stgf{R}$.
\end{chunk}

As $\gfc{R}$ is a Frobenius category, we may relate extensions in $\gfc{R}$ with morphisms in $\stgf{R}$: recall that there is a standard identification $\underline{\Hom}(X,Y)\simeq \t{Ext}_{R}^{1}(X,\Omega Y)$, where $\Omega$ is the inverse to the shift in $\stgf{R}$, see \cite[Lemma 3.3.3]{krbook}, for example. Given a map $\underline{f}\colon X\to \Sigma Y$ in $\stgf{R}$, let $\eta(\underline{f})$ denote the corresponding short exact sequence in $\msf{GFlat}(R)$ under this identification, which is obtained via considering the extension $0\to \Omega Y\to M\to X\to 0$ which arises from the triangle $\Omega Y\to M\to X\xrightarrow{f} Y$.

\begin{lem}(Modification of \cite[Proposition 4.2.2]{bengna})\label{pureexactrelationship}
Let $R$ be right coherent such that $\stgf{R}$ is compactly generated. Then the following are equivalent for $\underline{f}\colon M\to  N$ in $\stgf{R}$:
\begin{enumerate}
\item $\underline{f}$ is a phantom map;
\item $\eta(\underline{f})$ is pure relative to $\mathscr{C}$, that is $\Hom(C,\eta(\underline{f}))$ is exact for all $C\in\mathscr{C}$.
\end{enumerate}
\end{lem}

\begin{proof}
We modify the proof of \cite[Proposition 4.2.2]{bengna}. Suppose $\underline{f}$ is a map, which extends to a triangle
\[
\Omega N \overset{\underline{l}}{\longrightarrow} X\overset{\underline{m}}{\longrightarrow} M \overset{\underline{f}}{\longrightarrow}N 
\] 
that corresponds to an extension $0\to\Omega N \xrightarrow{l} X\xrightarrow{m} M\to 0$. 

Let us assume that $\underline{f}$ is phantom, then, for any compact object $C$, the induced map $\underline{\Hom}(C,\underline{f})$ is zero, hence $\underline{\Hom}(C,\underline{m})\colon\underline{\Hom}(C,X)\to\underline{\Hom}(C,M)$ is epic. We show that $\Hom(C,m)\colon\Hom(C,X)\to\Hom(C,M)$ is also epic, and thus the extension is pure relative to $\mathscr{C}$.

Now, suppose that $\alpha\in\Hom(C,M)$ and consider its image $\underline{\alpha}\in\underline{\Hom}(C,M)$, which is of the form $\underline{\alpha}=\underline{m}\underline{\beta}$ for some $\underline{\beta}\in\underline{\Hom}(C,X)$, the image of a map $\beta\colon C\to X$ in $\gfc{R}$. Consequently $\alpha-m\beta$ factors through a flat cotorsion module, say $F$. As $F$ is a projective object in $\gfc{R}$, it follows that the map $r\colon F\to M$, appearing in the factorisation of $\alpha-m\beta$, lifts along the surjection $m\colon X\to M$. In particular, we have $r=mp$ for some $p\colon F\to X$, and thus $\alpha-m\beta = mpq$ for some $q\colon C\to F$. In particular, $\alpha =m(\beta+pq)$ is in the image of $\Hom(C,m)$, which is what we wanted.

For the other implication, apply \cite[Lemma 3.3.3]{krbook} to see that if the extension $\eta(\underline{f})$ is pure with respect to $\mathscr{C}$, then $\underline{f}$ is immediately phantom.
\end{proof}

In particular, with this setup, it becomes apparent that the pure triangles in $\stgf{R}$ biject, through $\eta$, with the short exact sequences of Gorenstein flat $R$ modules which are pure relative to $\mathscr{C}$. 

In certain situations, the above lemma enables a direct comparison between the pure triangles in $\msf{GFlat}(R)$ which are termwise cotorsion, and the pure triangles in $\stgf{R}$.

\begin{prop}\label{extensionpureiff}
Let $R$ be a coherent ring that is cotorsion over itself and satisfies the following conditions:
\begin{enumerate}
\item $\msf{Gproj}(R)\subseteq\msf{GFlat}(R)$; hence $\msf{Gproj}(R)\subseteq\gfc{R}$;
\item $\msf{GFlat}(R)=\rlim\msf{Gproj}(R)$;
\item each finitely presented Gorenstein projective $R$-module is compact in $\stgf{R}$.
\end{enumerate}
Then $\stgf{R}$ is compactly generated by $\stgp{R}$. Furthermore, a triangle $X\to Y\to Z\to \Sigma X$ is pure in $\stgf{R}$ if and only if the corresponding exact sequence $0\to X\to Y\to Z\to 0$ is pure in $\msf{GFlat}(R)$.
\end{prop}

\begin{proof}
To begin with, observe that the assumption on the ring ensures that any finitely presented Gorenstein projective module $M$ has a complete acyclic resolution by finitely presented projective $R$-modules, see \cite[p. 257]{rha}, each of which, being a summand of a free module of finite rank, is cotorsion. Consequently, $M$ is a cycle in an acyclic complex of cotorsion modules, hence is itself cotorsion by \cite[Theorem 1.3]{bce}, which gives an inclusion $\msf{Gproj}(R)\subseteq\gfc{R}$, this justifies the hence in the first enumerated item. As such, there is an inclusion $\stgp{R}\subseteq\stgf{R}$, which illustrates $\stgp{R}$ as a triangulated subcategory of $\stgf{R}$.

Consider the colocalising subcategory 
\[
\D:=\{X\in\stgf{R}:\underline{\Hom}(C,X)=0 \t{ for all }C\in\stgp{R}\},
\]
and suppose $X\in\D$ is a pure injective object, which means, as $\stgp{R}$ is a triangulated subcategory, that $\Sigma^{i}X\in\D$ for all $i\in\Z$ as well. There is then a canonical isomorphism
\[
0\simeq \underline{\Hom}(M, X)\simeq \t{Ext}^{1}(M,\Sigma^{-1}X)\simeq \t{Ext}_{R}^{1}(M,\Sigma^{-1}X)
\]
for all $M\in\msf{Gproj}(R)$. Here, $\t{Ext}^{1}(M,X)$ is the group of Yoneda extensions in the exact structure on  $\gfc{R}$. Since $M$ is finitely presented and $\mod{R}\cap\msf{FlatCot}(R)=\msf{proj}(R)$, one obtains the isomorphism in the equation. We note that for an arbitrary object $Z\in\stgf{R}$, there is an inclusion of groups $\t{Ext}^{1}(Z,\Sigma^{-1}X)\hookrightarrow \t{Ext}_{R}^{1}(Z,\Sigma^{-1}X)$, given by seeing any extension in $\gfc{R}$ as an extension in $\Mod{R}$.

Now, the assumption that $\msf{GFlat}(R)=\rlim\msf{Gproj}(R)$ means we may express $X\simeq\rlim_{I}L_{i}$, with $L_{i}\in\msf{Gproj}(R)$. Then, as $X$, and so $\Sigma^{-1}X$, is a pure injective object in both $\stgf{R}$ and $\Mod{R}$ by \cref{piinst}, there are isomorphisms
\[
\underline{\Hom}(X,X)\simeq \t{Ext}^{1}(X,\Sigma^{-1}X)\subseteq \t{Ext}_{R}^{1}(X,\Sigma^{-1}X)\simeq \t{Ext}_{R}^{1}(\rlim_{I}L_{i},\Sigma^{-1}X)\simeq \llim_{I}\t{Ext}_{R}^{1}(L_{i},\Sigma^{-1}X)=0
\]
where the inclusion follows from the previous discussion, and the final isomorphism follows from \cite[Lemma 6.28]{gt}. In particular, it follows that if $X\in\D$ is pure injective, then $X=0$, or equivalently $X$ is flat when viewed as an $R$-module.

Now, suppose $N\in\D$, so we again have $\t{Ext}_{R}^{1}(M,N)=0$ for all $M\in\msf{Gproj}(R)$. As each such $M$ has an infinite left resolution by finitely presented projective $R$-modules, it follows that 
\[
\mc{X}:=\{Z\in\Mod{R}:\t{Ext}_{R}^{1}(M,Z)=0\t{ for all }M\in\msf{Gproj}(R)\}
\]
is a definable subcategory of $\Mod{R}$ by \cite[Theorem 10.2.35]{psl}, and thus $\mc{D}:=\mc{X}\cap\msf{GFlat}(R)$ is a definable subcategory of $\msf{GFlat}(R)$. In particular, if $N\in\D$ then, viewed as an $R$-module, $N\in\mc{D}$, and as $\mc{D}$ is closed under pure injective envelopes, we have that $PE(N)\in\mc{D}$, hence $PE(N)\in\D$ is pure injective when viewed as an object of $\stgf{R}$, again by \cref{piinst}. In particular, from the above we deduce that $PE(N)$ is a flat $R$-module, and thus $N$ is a flat $R$-module, as $R$ is coherent. In other words, we have shown $N\in\D$ if and only if $N=0$, hence $\stgp{R}$ generates $\stgf{R}$, and by the third assumption, these form a small triangulated subcategory of compact generators. As compact objects in enhanced triangulated categories are unique, it follows that $\stgf{R}^{\c}=\stgp{R}$.

Let us now move to the furthermore statements. Suppose that $X\to Y\to Z\to \Sigma X$ is a pure triangle so $Z\to \Sigma X$ is phantom. Then, by \cref{pureexactrelationship}, the extension $0\to X\to Y\to Z\to 0$ is pure with respect to the stable objects in $\msf{Gproj}(R)$. The other objects in $\msf{Gproj}(R)$ are projective, and therefore it follows that the extension is pure with respect to all of $\msf{Gproj}(R)$. This means, from the assumption that $\msf{GFlat}(R)=\rlim\msf{Gproj}(R)$, that the extension is pure in $\msf{GFlat}(R)$. 

For the converse, suppose that $X\to Y\to Z\xrightarrow{\underline{f}}\Sigma X$ is a triangle in $\stgf{R}$, such that the corresponding short exact sequence $\eta(\underline{f})=(0\to X\to Y\to Z\to 0)$ is pure in $\msf{GFlat}(R)$. By \cref{pureexactrelationship}, for $\underline{f}$ to be phantom it suffices to show that $\Hom_{R}(C,\eta(\underline{f}))$ is exact for each object that is compact in $\stgf{R}$. We showed earlier in the proof that $\stgf{R}^{\c}=\stgp{R}$. Yet the assumption that $\rlim\msf{Gproj}(R)=\msf{GFlat}(R)$ tells us that a short exact sequence $\mbb{S}$ is pure if and only if $\Hom_{R}(M,\mbb{S})$ is exact for all $M\in\msf{Gproj}(R)$, see \cite{CrawleyBoevey}. In particular, the assumption that $\eta(f)$ is pure, tells us that $\Hom_{R}(\msf{Gproj}(R),\eta(f))$ is exact, hence $f$ is phantom and the triangle is pure.
\end{proof}

\begin{rem}
Rings where the assumptions of \cref{extensionpureiff} hold appear quite naturally. For example, any Iwanaga-Gorenstein ring $R$ when $R$ is cotorsion, for instance a complete commutative Gorenstein ring, or a Gorenstein Artin algebra. Another example comes from local Cohen-Macaulay local rings admitting a dualising module, see \cite[Theorem 10.4.31]{rha}.
\end{rem}

Let us now assume that $R$ is a ring which satisfies the conditions of \cref{extensionpureiff}. Given a definable subcategory $\mc{D}\subset\msf{GFlat}(R)$, let $\underline{\mc{D}}$ denote its image in $\stgf{R}$; in other words, the non-zero objects in $\underline{\mc{D}}$ are just the non-flat cotorsion objects in $\mc{D}$.

\begin{lem}\label{stdefinable}
With the assumptions of \cref{extensionpureiff}, if $\mc{D}\subseteq\msf{GFlat}(R)$ is definable, then $\underline{\mc{D}}\subseteq\stgf{R}$ is definable.
\end{lem}
\begin{proof}
Since $\stgf{R}$ is an algebraic triangulated category, we may, by \cite[Theorem 4.7]{lv}, show that $\underline{\mc{D}}$ is closed under pure subobjects, pure quotients and products to deduce it is definable. As products in $\stgf{R}$ are identical to those in $\msf{GFlat}(R)$, and the fact that both $\mc{D}$ and $\msf{Cot}(R)$ are closed under products, it follows immediately that $\underline{\mc{D}}$ is closed under products; note the product of non-flat objects is non-flat, as each object is a summand, and thus pure, in the product. 

We now show closure under pure subobjects and pure quotients. To this end, suppose $X\to Y\to Z\to \Sigma X$ is a pure triangle in $\stgf{R}$ with $Y\in\underline{\mc{D}}$. Then, by \cref{extensionpureiff} and the assumptions on $R$, it follows that $0\to X\to Y\to Z\to 0$ is a pure exact sequence in $\msf{GFlat}(R)$. In particular, $X$ and $Z$ must lie in $\mc{D}$ as $\mc{D}$ is closed under pure subobjects and quotients, which is what we wanted.
\end{proof}

Let us now compare the Ziegler spectra of $\stgf{R}$ and $\msf{GFlat}(R)$.

\begin{thm}\label{zieglerdisjoint}
Let $R$ satisfy the conditions of \cref{extensionpureiff}, then there is a homeomorphism
\[
\msf{Zg}(\msf{GFlat}(R))\simeq \msf{Zg}(\stgf{R})\coprod \msf{Zg}(\Flat{R}).
\]
\end{thm}

\begin{proof}
By \cref{piinst}, an object of $\stgf{R}$ is indecomposable pure injective if and only if it is indecomposable and pure injective in $\msf{GFlat}(R)$. In particular, the identity map gives an inclusion $i\colon \msf{pinj}(\stgf{R})\to\msf{pinj}(\msf{GFlat}(R))$. This map is continuous: indeed, if $X\subset\msf{Zg}(\msf{GFlat}(R))$ is closed, then, by the definition of the Ziegler topology, there is a definable subcategory $\mc{D}\subseteq\msf{GFlat}(R)$ such that $X=\mc{D}\cap\msf{pinj}(\msf{GFlat}(R))$. By \cref{stdefinable}, the class $\underline{\mc{D}}$ is a definable subcategory, and the indecomposable pure injective objects in $\underline{\mc{D}}$ are precisely $i^{-1}(X)$, hence the map is continuous as claimed.

Meanwhile, the inclusion $j\colon\Flat{R}\to\msf{GFlat}(R)$ clearly induces a continuous map $\msf{Zg}(\Flat{R})\to\msf{Zg}(\msf{GFlat}(R))$ by \cite[Theorem 15.5]{dac}.

These two continuous maps give us a map $\msf{Zg}(\msf{Flat}(R))\sqcup\msf{Zg}(\stgf{R})\to\msf{Zg}(\msf{GFlat}(R))$. The statement of the theorem holds provided there is a homeomorphism $\msf{Zg}(\stgf{R})\simeq \msf{Zg}(\msf{GFlat}(R))\setminus\msf{Zg}(\Flat{R})$. Yet this holds: with the assumptions of \cref{extensionpureiff}, the category $\stgf{R}$ is compactly generated, with compact objects $\stgp{R}$. By \cite[Lemma 2.9]{birdwilliamsonhomological}, there is a homeomorphism $\msf{Zg}(\stgf{R})\simeq \msf{Zg}(\rlim\stgp{R})$, where $\rlim\stgp{R}$ is viewed as the category of ind-objects over $\stgp{R}$, or equivalently flat (i.e. cohomological) functors $\stgp{R}^{\t{op}}\to\ab$. Yet $\msf{Zg}(\rlim\stgp{R})\simeq\msf{Zg}(\msf{GFlat}(R))\setminus\msf{Zg}(\Flat{R})$ by the results of \cite[\S 5]{kredc}. 
\end{proof}

We shall use this result in \cref{section:knorrer}, when periodicity for Gorenstein flat modules over hypersurfaces is considered.

\begin{ex}
Let us describe how we can see $\stgf{R}$ as a measure of the discrepancy of the inclusions $\msf{proj}(R)\subseteq\msf{Gproj}(R)$ and $\msf{Flat}(R)\subseteq\msf{GFlat}(R)$. This shows that, from a perspective of purity, $\stgf{R}$ is the suitable unbounded analogue of $\stgp{R}$, instead of, for example, $\ul{\msf{GProj}}(R)$.

Over any ring, there is an inclusion of classes $\msf{PGFlat}(R)\subseteq\msf{GProj}(R)\subseteq\msf{GFlat}(R)$, where $\msf{PGFlat}(R)$ are the projectively coresolved Gorenstein flat modules introduced in \cite{sarstov}, which we do not introduce in more detail since they are not used elsewhere. In particular, there is then an inclusion of definable closures $\msf{Def}(\msf{PGFlat}(R))\subseteq\msf{Def}(\msf{GProj}(R))\subseteq\msf{Def}(\msf{GFlat}(R))$.

Now, it is stated on \cite[p. 28]{sarstov} that $\msf{GFlat}(R)$ is closed under pure quotients if and only if the closure of $\msf{PGFlat}(R)$ under pure quotients is $\msf{GFlat}(R)$. As such, whenever $\msf{GFlat}(R)$ is definable, we have that $\msf{GFlat}(R)\subseteq \msf{Def}(\msf{PGFlat}(R))$, hence the two are equal, and thus also equal to the definable closure of $\msf{GProj}(R)$. Of course, the assumption that $\msf{GFlat}(R)$ is definable forces the ring to be right coherent and for $\msf{Flat}(R)$ to be definable.

Now, assuming that the definable closure of $\msf{Gproj}(R)$ is $\msf{GFlat}(R)$, which holds whenever the assumptions of \cref{extensionpureiff} do, we have an inclusion $\msf{Flat}(R)\subseteq\msf{GFlat}(R)$ of definable categories, build from definably closing the inclusion $\msf{proj}(R)\subseteq\msf{Gproj}(R)$.

We may then consider the \ti{definable quotient} category corresponding to the inclusion $\msf{Flat}(R)\subseteq\msf{GFlat}(R)$. Let us briefly recall the construction of this category. There are Serre subcategories $\mc{S}_{\msf{F}}$ and $\mc{S}_{\msf{GF}}$ of $\msf{fp}(\mod{R},\ab)$ such that $\Flat{R}=\{X\in\Mod{R}:\bar{f}X=0\t{ for all }f\in\mc{S}_{\msf{F}}\}$, and likewise for $\msf{GFlat}(R)$. As $\Flat{R}\subset\msf{GFlat}(R)$, there is a reverse inclusion $\mc{S}_{\msf{GF}}\subset\mc{S}_{\msf{F}}$ of small abelian categories, and by localising there is a resulting small abelian category $\mc{S}_{\msf{GF/F}}$. It follows, from \cite{kredc} and \cite[Chapter 11]{dac}, that there is a unique definable category $\mc{D}$, which is the definable quotient category corresponding to the inclusion $\msf{Flat}(R)\subseteq\msf{GFlat}(R)$.

While $\mc{D}$ is a definable subcategory, it is not a definable subcategory of $\Mod{R}$. However, it has the property that 
\[
\msf{Zg}(\mc{D})\simeq \msf{Zg}(\msf{GFlat}(R))\setminus \msf{Zg}(\Flat{R})
\]
by \cite[Corollary 6.3]{kredc}. By \cref{piinst}, we see that the points of $\msf{Zg}(\mc{D})$ are precisely the points of $\msf{Zg}(\stgf{R})$. When the conditions of \cref{extensionpureiff} hold, we see that $\msf{Zg}(\mc{D})$ is actually homeomorphic to $\msf{Zg}(\stgf{R})$. 

In this case, we can actually explicitly give $\mc{D}$ as the finitely accessible category $\msf{Flat}(\stgp{R}^{\op})$, which is just the ind completion of $\stgp{R}$, by applying \cite[Theorem 5.4]{kredc}
\end{ex}

\begin{ex}
We now give an explicit example of the Ziegler spectrum of $\stgf{R}$ over a hypersurface singularity. Let $k$ be an algebraically closed field of characteristic other than 2, and consider the $A_{\infty}$ curve singularity $R=k[[x,y]]/x^2$.

These rings is of interest in the study of maximal Cohen-Macaulay modules, as it was shown in \cite[Theorem B]{bgs} that a one-dimensional curve singularity admits countably many indecomposable maximal Cohen-Macaulay modules if and only if it is either the $A_{\infty}$, which we consider, or the $D_{\infty}$ curve singularity. 

As $R$ is Gorenstein, we may consider the category $\msf{K}_{\t{ac}}(\msf{Inj}(R))$. There is a chain of triangulated equivalences 
\[
\stgf{R}\simeq\msf{K}_{\t{tac}}(\msf{Proj}(R))\simeq\msf{K}_{\t{ac}}(\msf{Inj}(R)),
\]
where the categories are all compactly generated, and the compact objects is the singularity category $\msf{D}_{\t{sg}}(R)$. By the preceding discussions, the indecomposable pure injective objects of this category are the non-flat indecomposable pure injective objects in $\msf{GFlat}(R)$.

The finitely generated indecomposable pure injective Gorenstein flat modules over $R$ are just the indecomposable maximal Cohen-Macaulay modules, which were classified at at \cite[\S4]{bgs}. These are the ideals $I_{n}:=(x,y^{n})$, with $n\geq 0$ and $I_{\infty}=xR$. Clearly $R=I_{0}$ is the only flat-cotorsion finitely generated maximal Cohen-Macaulay module.

By \cite[Theorem 3.8]{bird}, it can be deduced that the non-finitely generated indecomposable pure injectives were described at \cite[Theorem 6.1]{puninski}, where there is no mention that the described modules are, in fact, Gorenstein flat. The non-finitely generated objects are $\bar{R}$, the normalisation of $R$, $Q(R)$, the total ring of quotients for $R$, and the Laurent series ring $k((y))$ viewed as an $R$-module along the projection $R\to k[[y]]$, where $k[[y]]\simeq R/(x)$. Clearly $Q(R)$ is flat as it is a localisation, while $\bar{R}$ is most certainly not flat as it is a normalisation of a non-reduced ring. The module $k((y))$ is also not flat: let $\mf{p}=(x)$ be the prime ideal generated by $x$, and consider the indecomposable maximal Cohen-Macaulay module $I_{\infty}=xR$. Then $(I_{\infty})_{\mf{p}}\simeq k((y))$, which is not free over $R_{\mf{p}}$. In particular, $k((y))=k((y))_{\mf{p}}$ is not flat over $R_{\mf{p}}$, and as flatness is a local property, it follows that $k((y))$ cannot be flat over $R$. 

Therefore the indecomposable pure injective objects in $\stgf{R}$ are the non-free indecomposable maximal Cohen-Macaulay modules $I_{n}$, for $1\leq n\leq \infty$, $\bar{R}$ and $k((y))$.

In \cite[Proposition 6.4]{puninski}, it was shown that the isolated points of $\msf{Zg}(\msf{GFlat}(R))$ are the $I_{n}$, for $n<\infty$. By \cref{stdefinable}, it follows that $\underline{I}_{n}$, the image of $I_{n}$ in $\stgf{R}$, is also isolated, for $0<n<\infty$.

Now, by \cite[Theorem 6.9]{puninski} and the comments following, $k((y))$ is a closed point of $\msf{Zg}(\msf{GFlat}(R))$, and thus, again by \cref{stdefinable}, $\underline{k((y))}$ is again closed. The remaining points, $\underline{I}_{\infty}$ and $\underline{\bar{R}}$ are neither closed, nor isolated - they become isolated once the points $\underline{I}_{n}$ are removed, by \cite[Lemma 6.7]{puninski}.
\end{ex}

\begin{chunk}
Let us finish this section by showing that whenever $R$ is right coherent such that every flat module has finite projective dimension every pure triangle in $\stgf{R}$ arises from a triangle in $\msf{K}(\msf{Proj}(R))$, which is compactly generated whenever $R$ is right coherent by \cite{neemanhomotopy}. 

For such rings, there is a sequence of triangulated equivalences of categories
\[
\begin{tikzcd}[column sep=1.5cm]
\msf{K}_{\t{tac}}(\msf{Proj}(R)) \arrow[r,,draw=none,"\scriptstyle\perp" marking] \arrow[r, shift left = 1ex, "Z_{0}"] \arrow[r, leftarrow, shift right = 1ex, "P", swap] & \stgP{R}  \arrow[r,,draw=none,"\scriptstyle\perp" marking] \arrow[r, shift left = 1ex, "\alpha"] \arrow[r, leftarrow, shift right = 1ex, "\beta", swap]& \stgf{R}  \arrow[r, shift right = 1ex, "F", , swap] \arrow[r, leftarrow, shift left = 1ex, "Z_{0}"]\arrow[r,,draw=none,"\scriptstyle\perp" marking]&\msf{K}_{\t{tac}}(\msf{FlatCot}(R)) 
\end{tikzcd}
\]
where $\alpha=\msf{Ho}(\t{Id})$, by \cite[Theorem 3.8]{CET} and \cite[Theorem 5.1]{estgil}.

Suppose $R$ is Noetherian with a dualising complex. Then there is, by \cite[Theorem 5.6.1]{krloc} and \cite[Theorem 5.3]{ik} (see also \cite[Proposition 2.1.5]{ops}), a recollement of compactly generated triangulated categories
\[
\begin{tikzcd}[column sep=2cm]
\msf{Loc}(R, \rho\Hom_{R}(\omega_{R},\lambda R)) \arrow[r, shift left = 2ex, "I=\mrm{inc}"] \arrow[r, leftarrow, "I^{*}" description] \arrow[r, shift right = 2ex, "I_{*}", swap] & \msf{K}(\msf{Proj}(R))  \arrow[r, shift left = 2ex, "Q"] \arrow[r, leftarrow, "Q^{*}" description] \arrow[r, shift right = 2ex, "Q_{*}", swap]  & \msf{K}_{\t{tac}}(\msf{Proj}(R))
\end{tikzcd}
\]
where $\rho\colon\msf{K}(R)\to\msf{K}(\msf{Proj}(R))$ is the right adjoint to the coproduct preserving inclusion, and $\lambda\colon\D(R)\to\msf{K}(R)$ is an injective resolution. In particular, $Q$ endows $\msf{K}_{\t{tac}}(\msf{Proj}(R))$ with a set of compact generators.

Now, if $X\to Y\to Z\to \Sigma X$ is a pure triangle in $\stgf{R}$, then, as equivalences preserve purity, the triangle $P\beta X\to P\beta Y\to P\beta Z\to \Sigma P\beta X$ is pure in $\msf{K}_{\t{tac}}(\msf{Proj}(R))$. As $Q^{*}$ has a left and right adjoint, it preserves all homotopy colimits, which exist in $\msf{K}_{\t{tac}}(\msf{Proj}(R))$ as it has an $\infty$-categorical enhancement. In particular, $Q^{*}$ preserves pure triangles. Thus 
\[
Q^{*}P\beta X\to Q^{*}P\beta Y\to Q^{*}P\beta Z\to \Sigma Q^{*} P\beta X
\]
is a pure triangle in $\msf{K}(\msf{Proj}(R))$. Applying $\alpha Z_{0}Q$ to this triangle, and observing that $QQ^{*}=\t{Id}$, we see that we recover the original triangle. This shows that every pure triangle in $\stgf{R}$ is the image of a pure triangle in $\msf{K}(\msf{Proj}(R))$. 

\end{chunk}

\section{Ascent for Gorenstein flat cotorsion modules}\label{ascentsection}

Let us suppose that $R\to S$ is a ring map, such that both $R$ and $S$ are right coherent and $S$ is finitely presented as a right $R$-module. Then the extension of scalars functor $S\otimes_{R}-\colon\Mod{R}\to\Mod{S}$ is a definable functor by \cref{extensionflatexample}. Moreover, this functor preserves flat modules, so gives a definable functor $\Flat{R}\to\Flat{S}$. 

If we additionally assume that $S$ has finite flat dimension over $R$, then $S\otimes_{R}-$ preserves Gorenstein flat modules by \cite[5.6(a)]{ascent}. In particular, if we assume that $\msf{GFlat}(R)$ and $\msf{GFlat}(S)$ are definable, we also have a definable functor $\msf{GFlat}(R)\to\msf{GFlat}(S)$.

In this section, we investigate what this setup yields on the level of stable categories, and its implications for purity.

The first goal of this section is to show that this induces a triangulated functor $\stgf{R}\to\stgf{S}$. The following lemma lays the foundation for such a statement.

\begin{lem}\label{gflatextscal}
Let $\phi\colon R\to S$ be a morphism of right coherent rings such that $S\in\mod{R^{\op}}$ and is of finite flat dimension over $R$. Then extension of scalars along $\phi$ induces a functor
\[
S\otimes_{R}-\colon\msf{Ch}_{\mrm{tac}}(\msf{FlatCot}(R))\to\msf{Ch}_{\mrm{tac}}(\msf{FlatCot}(S)).
\]
In particular, we obtain a functor $S\otimes_{R}-\colon\gfc{R}\to\gfc{S}$.
\end{lem}

\begin{proof}
Let $T$ be a totally acyclic complex of flat cotorsion $R$-modules, and consider the complex $S\otimes_{R}T\in\msf{Ch}(S)$; we show that the latter satisfies the conditions of being a totally acyclic complex of flat cotorsion $S$-modules, as detailed in \cref{totallyacyclicdef} for $\msf{U}=\msf{FlatCot}(S)$. 

Since, as stated above, $S\otimes_{R}-\colon\Flat{R}\to\Flat{S}$ is a definable functor with the assumptions on $S$, it preserves pure injective, and thus flat cotorsion modules by \cref{Flatcotorsionpair}. In particular, the second condition of \cref{totallyacyclicdef} is obvious.

What remains is to prove the acyclicity conditions. Let us first show that $S\otimes_{R}T$ is itself acyclic. By \cite[Lemma A.7]{dcmca}, we must show that $\t{Tor}_{i}^{R}(S,T_{j})=0$ for all $i>0$ and $j\in\Z$ and $\t{Tor}_{i}^{R}(S,\msf{Z}_{j}(T))=0$ for all $i>0$ and $j\ll 0$. The first is trivial, since $T_{j}\in\Flat{R}$ by definition. For the second Tor vanishing, for all $i>0$ and $j\in\Z$ as $\msf{Z}_{j}(T)$ is a Gorenstein flat module by \cref{gorensteinflat}, and $S$ being finite flat dimension means $\t{Tor}_{i}^{R}(S,\msf{Z}_{j}(T))=0$ for all $j<0$ by \cite[Lemma 2.3]{cfh}.

Having just shown that $S\otimes_{R}T$ is an acyclic complex of flat cotorsion $S$-modules, \cite[Theorem 1.3]{bce} tells us that $\msf{Z}_{j}(S\otimes_{R}T)\in\msf{Cot}(S)$ for all $j\in\Z$. Consequently, as $\msf{FlatCot}(S)$ is self orthogonal, we have that $\t{Ext}_{S}^{i}(F,(S\otimes_{R}T)_{j})=0=\t{Ext}_{S}^{i}(F,\msf{Z}_{j}(S\otimes_{R}T))$ for all $i>0$ and $j\in\Z$, and $F\in\msf{FlatCot}(S)$; consequently, by \cite[Lemma A.4]{dcmca}, $\Hom_{S}(F,S\otimes_{R}T)$ is acyclic.

Lastly, we show $\Hom_{S}(S\otimes_{R}T,F)$ is acyclic for all $F\in\msf{FlatCot}(S)$. Since $S\otimes_{R}-\colon\msf{Ch}(R)\to\msf{Ch}(S)$ arises as the extension of a right exact functor $\Mod{R}\to\Mod{S}$, there are isomorphisms $\msf{C}_{j}(S\otimes_{R}T)\simeq S\otimes_{R}\msf{Z}_{-j}(T)$ for all $j\in\Z$ by \cite[E 2.2.7]{dcmca}. By the assumption of $R$ being right coherent, we have that $\msf{Z}_{j}(T)$ is Gorenstein flat for all $j\in\Z$, and thus $S\otimes_{R}\msf{Z}_{j}(T)\in\msf{GFlat}(S)$ by \cite[Lemma 4.6.(i)]{ascent}.  Recall from \cref{GFlatcotorsionpair} that $\msf{GFlat}(R)\cap\msf{GCot}(R)=\msf{FlatCot}(R)$, and therefore $\t{Ext}_{S}^{i}(S\otimes_{R}\msf{Z}_{j}(T),F)=0$ for all $i>0$ and $j\in\Z$, and $F\in\msf{FlatCot}(R)$. Moreover, for the same reason, $\t{Ext}_{S}^{i}((S\otimes_{R}T)_{j},F)=0$. This concludes the proof.
\end{proof}

\begin{rem}
The assumptions of the above lemma can be relaxed to some extent. One can replace $S$ with any $S-R$-bimodule $M$ such that $M$ is quasi-isomorphic in $\D(R)$ to a perfect complex, and is flat over $S$. In such a situation, $M\otimes_{R}-\colon\Flat{R}\to\Flat{S}$ is still a definable functor that preserves Gorenstein flat modules, by \cite[Lemma 4.6]{ascent}.
\end{rem}

\begin{ex}\label{regseq}
Let us give an example of a ring homomorphism which satisfies the conditions of \cref{gflatextscal} which shall appear throughout the rest of the article. Suppose that $R$ is a commutative Noetherian ring, and consider a sequence of elements $\bm{x}=x_{1},\cdots,x_{n}$ of $R$. Provided the following two conditions hold, the sequence $\bm{x}$ is called an $R$-\ti{sequence}:
\begin{enumerate}
\item For each $0<i<n$ the multiplication map $\cdot x_{i+1}\colon R/(x_{1},\cdots,x_{i})R\to R/(x_{1},\cdots,x_{i})R$ is injective;
\item $R/\bm{x}R\neq 0$.
\end{enumerate}
For instance, in $R=k[x_{1},\cdots,x_{n}]$, where $k$ is a domain, any subset of the indeterminates forms an $R$-sequence. 

Given an $R$-sequence $\bm{x}$, we can consider the canonical projection map $R\to R/\bm{x}$. As $R$ is Noetherian, it is clear that $R/\bm{x}$ is finitely presented. To show that $R/\bm{x}$ has finite flat dimension, we note that, in $\D(R)$, there is a quasi-isomophism between $R/\bm{x}$ and the \ti{Koszul complex} on $\bm{x}$. 

The Koszul complex on a single element $r\in R$ is the two-term complex $K(r,R)=0\to R\xrightarrow{\cdot r}R\to 0$, and the Koszul complex on an arbitrary sequence $\bm{r}=r_{1},\cdots,r_{m}$ is given by $K(\bm{r},R)=K(r_{1},R)\otimes K(r_{2},R)\otimes\cdots\otimes K(r_{m},R)$. Clearly, $K(\bm{r},R)$ is a perfect complex. And the quotient ring $R/\bm{r}$ is quasi-isomorphic in $\D(R)$ to $K(\bm{r},R)$ if and only if $\bm{r}$ is an $R$-sequence.
\end{ex}

Let us now show that the conclusions of \cref{gflatextscal} can be extended to functors between triangulated categories.

\begin{thm}\label{triangulatedfunctor}
Let $R\to S$ be a morphism of right coherent rings such that $S\in\mod{R^{\op}}$ and is of finite flat dimension over $R$. Then $S\otimes_{R}-\colon\gfc{R}\to\gfc{S}$ is an exact functor of Frobenius categories, and the induced functor
\[
S\otimes_{R}-\colon\stgf{R}\to\stgf{S}
\]
is triangulated.
\end{thm} 

\begin{proof}
The existence of the functor $S\otimes_{R}-\colon\gfc{R}\to\gfc{S}$ was proved in \cref{gflatextscal}. The assumptions on $R\to S$ ensure this restricts to a functor $\msf{FlatCot}(R)\to\msf{FlatCot}(S)$; these are precisely the projective-injective objects in these Frobenius categories. Moreover, the functor is exact, for $\t{Tor}_{1}^{R}(S,X)=0$ for any Gorenstein flat $R$-module $X$, by the assumption that $S$ is of finite flat dimension and \cite[Lemma 4.6.(1)]{ascent}. 

Consequently, there is, by \cite[\S2.8]{Happel}, an induced functor $\stgf{R}\to\stgf{S}$, which is additive. The fact that it is triangulated follows from the fact that $S\otimes_{R}-\colon\msf{K}(\msf{FlatCot}(R))\to\msf{K}(\msf{FlatCot}(S))$ is triangulated, and combining the above with the triangulated equivalence
\[
\stgf{R}\simeq\msf{K}_{\t{tac}}(\msf{FlatCot}(R))
\]
shows that the restriction to  $\stgf{R}\to\stgf{S}$ is also triangulated.
\end{proof}

We now show that the induced functor on stable categories preserves coproducts and products.

\begin{lem}\label{preservesprodsandcoprods}
Let $R\to S$ be a morphism of coherent rings such that $S\in\mod{R^{\op}}$ and is of finite flat dimension over $R$. Then the functor $S\otimes_{R}-\colon\stgf{R}\to\stgf{S}$ of \cref{triangulatedfunctor} preserves products. Moreover, if every flat $R$-module has finite projective dimension, and the same for flat modules over $S$, the functor additionally preserves coproducts.
\end{lem}

\begin{proof}
Preservation of products is immediate from the assumption on $S$ as an $R$-module and the fact that $S\otimes_{R}-\colon\Mod{R}\to\Mod{S}$ preserves products. We now consider preservation of coproducts. By the assumptions of every flat $R$ and $S$-module having finite projective dimension, the discussion in \cref{gprojmodel} and \cite[Theorem 5.1]{estgil}, there is a diagram 
\[
\begin{tikzcd}
\Mod{R}_{\msf{GProj}} \arrow[r, "\t{Id}"] \arrow[d, "S\otimes_{R}-"]& \Mod{R}_{\msf{GFlat}}  \arrow[d, "S\otimes_{R}-"]\\
\Mod{S}_{\msf{GProj}} \arrow[r, "\t{Id}"]& \Mod{S}_{\msf{GFlat}}
\end{tikzcd}
\]
where the subscript indicates the model structure put on the category, which were described in \cref{coproducts} and \cref{gprojmodel}. It was also shown in \cite[Theorem 5.1]{estgil} that the horizontal arrows are (left) Quillen equivalences. We claim the vertical arrows are left Quillen functors. 

The right adjoint to $S\otimes_{R}-$ is $\msf{res}$, the restriction of scalars functor, realisable as $\Hom_{S}(S,-)\colon\Mod{S}\to\Mod{R}$. To show that the vertical arrows are left Quillen, we must show that $S\otimes_{R}-$ preserves cofibrations and $\msf{res}$ preserves fibrations.

For the Gorenstein projective model structure, the assumption that $S$ is finitely presented of finite flat dimension over $R$, which is right coherent, means that every $S$ module of finite projective dimension also has finite projective dimension over $R$. We may therefore apply \cite[5.6(b)]{ascent} to see that $S\otimes_{R}-$ preserves Gorenstein projective modules. The assumption that every flat $R$-module has finite projective dimension means that every Gorenstein projective module is Gorenstein flat, and thus $\t{Tor}_{i}^{R}(S,P)=0$ for each $P\in\msf{GProj}(R)$ by \cite[Lemma 2.3]{cfh}. In particular, $S\otimes_{R}-$ preserves cofibrations, which are monomorphisms with Gorenstein projective cokernel. A very similar argument shows that $S\otimes_{R}-$ preserves cofibrations in the $\msf{GFlat}$ model structure.

Since the fibrations in the $\msf{GProj}$ model structure on $\Mod{S}$ are all exact sequences, that $\msf{res}$ preserves fibrations is trivial, as it is exact. The fibrations in the $\msf{GFlat}$ model structure are epimorphisms with cotorsion kernel. It suffices to show that $\msf{res}$ preserves cotorsion modules, as it is exact. But it is clear that there are isomorphisms
\[
\t{Ext}_{R}^{i}(X,\msf{res}(Y)\simeq\t{Ext}_{S}^{i}(S\otimes_{R}X,Y)
\]
for all flat $R$-modules $X$ and all $S$-modules $Y$. In particular, if $X$ is flat and $Y$ is cotorsion, we have that $S\otimes_{R}X\in\Flat{S}$, and thus both Ext groups vanish for all $i>0$, hence $\msf{res}$ preserves cotorsion objects. This proves the claim that the vertical arrows are left Quillen.

Consequently, the induced square on homotopy categories
\[
\begin{tikzcd}
\ul{\msf{GProj}}(R) \arrow[r, "\simeq"] \arrow[d, "S\otimes_{R}-", swap]& \stgf{R}\arrow[d, "S\otimes_{R}-"] \\
\ul{\msf{GProj}}(S) \arrow[r, "\simeq"]& \stgf{S}
\end{tikzcd}
\]
commutes. In particular, we can deduce that $S\otimes_{R}-$ preserves coproducts, as the coproducts on $\stgf{R}$ are those induced from $\ul{\msf{GProj}}(R)$ and likewise for $S$.
\end{proof}

\begin{rem}
It follows that whenever $\stgf{R}$ and $\stgf{S}$ are compactly generated, the functor $S\otimes_{R}-\colon\stgf{R}\to\stgf{S}$ is a definable functor in the sense of \cite{definablefunctors}.
\end{rem}

\begin{chunk}\label{puritygobetween}
Let us consider how \cref{triangulatedfunctor} and \cref{preservesprodsandcoprods} enable us to discuss how purity is transported. Suppose that $R$ and $S$ are rings such that both $\msf{GFlat}(R)$ and $\msf{GFlat}(S)$ are definable subcategories. If $R\to S$ is a ring map satisfying the conditions of \cref{triangulatedfunctor}, there is a commuting diagram
\begin{equation}\label{diagram}
\begin{tikzcd}
\msf{GFlat}(R) \arrow[r] & \msf{GFlat}(S) \\
\gfc{R} \arrow[u, hook] \arrow[d, "\pi"] \arrow[r] & \gfc{S} \arrow[u, hook] \arrow[d, "\pi"] \\
\stgf{R} \arrow[r] & \stgf{S}
\end{tikzcd}
\end{equation}
where all three right arrows are $S\otimes_{R}-$, and $\pi$ is the canonical functor. In the top row, the functor is definable, and therefore, by \cref{definablefunctor}, it preserves the pure structure - namely the pure exact sequences and pure injective objects in $\msf{GFlat}(R)$. 

Now, by the discussion in \cref{Flatcotorsionpair} and \cref{piinst}, an object in $\msf{GFlat}(R)$ is pure injective if and only if it is pure injective in $\gfc{R}$ if and only if it is pure injective in $\stgf{R}$. Consequently, on the level of pure injective non-flat objects, there is no difference between any of the functors in the above diagram.

Thus, to understand the behaviour of the definable functor $\msf{GFlat}(R)\to\msf{GFlat}(S)$ on pure injective objects, it is enough to understand the behaviour of the functor $\stgf{R}\to\stgf{S}$. As illustrated in \cref{preservesprodsandcoprods}, the functor $S\otimes_{R}-$ preserves products, and therefore also preserves pure-injective objects. 

However, we would also like to understand whether the properties $S\otimes_{R}-\colon\stgf{R}\to\stgf{S}$ has, in relation to morphisms between pure injective objects, can be used to yield analogous information for the definable functor $S\otimes_{R}-\colon\msf{GFlat}(R)\to\msf{GFlat}(S)$.
\end{chunk}

We now turn our attention to conditions on $R\to S$ which ensure that this information can be lifted. Before stating the pertinent result, we require some auxiliary results.

\begin{lem}\label{psoofimage}
Suppose $R\to S$ be a map of right coherent rings such that $S$ is finitely presented as a right $R$-module. If $F\in\Flat{S}$, then $F$ is a pure submodule of $S\otimes_{R}X$, with $X\in\msf{FlatCot}(R)$. If $F\in\msf{FlatCot}(S)$, then this pure embedding is split.
\end{lem}

\begin{proof}
As $R$ is coherent, we may view $\Flat{R}=\msf{Def}(R)$, where $\msf{Def}(-)$ denotes the definable closure (which is closing a class under direct limits, products and pure subobjects), and we may do the same for $S$. Let $\mc{I}$ denote the essential image of $\Flat{R}$ under the definable functor $S\otimes_{R}-$. Note that, while $\mc{I}$ is closed under direct limits and products, as $S\otimes_{R}-$ preserves them, it is not itself definable, as it need not be closed under pure subobjects.

However, $\msf{pure}(\mc{I})$, the closure of $\mc{I}$ under pure subobjects is definable, see the remark after \cite[Corollary 13.4]{dac}. In particular, $\msf{pure}(\mc{I})=\Flat{S}$, as $\mc{I}\subset\Flat{S}$ and the latter is definable, and clearly $\Flat{S}\subseteq\msf{pure}(\mc{I})$ as the $S\in\mc{I}$. Therefore, for any flat module $F\in\Flat{S}$, there is a pure embedding $0\to F\to S\otimes_{R}X$ with $X\in\Flat{R}$.

We now show that $X$ can be chosen to be pure injective, and thus an object of $\msf{FlatCot}(R)$. To this end, consider the pure injective envelope $X\to PE(X)$ which is in $\Flat{R}$ by \cref{definableclosedunderpienvelopes}. Since $S\otimes_{R}-$ is definable, $S\otimes_{R}X\to S\otimes_{R}PE(X)$ is also a pure embedding, and thus $F\to S\otimes_{R}X\to S\otimes_{R}PE(X)$ is a composition of pure embeddings, which is pure, see \cite[Lemma 2.1.2]{psl}.

The second claim follows from the fact that if $F$ is in $\msf{FlatCot}(R)$ is it pure injective, and any pure monomorphism with a pure injective domain splits, see \cref{pureinjectivemodule}.
\end{proof}

The following corollary is now rather straightforward, and shows that, in the passage from $\gfc{S}$ to $\stgf{S}$, we need not check a morphism factors through any flat cotorsion $S$-module, just those in the image of $S\otimes_{R}-$. We assume the assumptions of \cref{psoofimage} are still in place.

\begin{cor}\label{factorsthroughimage}
Suppose $M$ and $N$ are $S$-modules. A morphism $f\colon M\to N$ factors through a flat cotorsion $S$-module if and only if it factors through $S\otimes_{R}X$ for some $X\in\msf{FlatCot}(R)$.
\end{cor}

\begin{proof}
Suppose $F\in\msf{FlatCot}(S)$, which is a summand of $S\otimes_{R}X$ for some $X\in\msf{FlatCot}(R)$ by \cref{psoofimage}. If $f$ factors through $F$, then it factors through $S\otimes_{R}X$ by passing through the splitting map $F\to S\otimes_{R}X$ and its retract. The other direction is trivial, since $S\otimes_{R}X$ is in $\msf{FlatCot}(S)$ for any $X\in\msf{FlatCot}(R)$.
\end{proof}

\begin{chunk}\label{fullconditions}
Consider the following to `fullness' conditions where $X\in\msf{FlatCot}(R)$ and $M\in\msf{GFlat}(R)\cap\msf{Pinj}(R)$:
\begin{enumerate}
\item[(C1)] the canonical map $\Hom_{R}(M,X)\to\Hom_{S}(S\otimes_{R}M,S\otimes_{R}X)$ is surjective;
\item[(C2)] the canonical map $\Hom_{R}(X,M)\to\Hom_{S}(S\otimes_{R}X,S\otimes_{R}M)$ is surjective.
\end{enumerate}
These conditions will be sufficient for us to lift properties from $\stgf{R}$ to $\msf{GFlat}(R)$. While, initially, the conditions appear quite strong, we shall show that in certain situations they hold for free, but before then, let us show what these conditions yield.
\end{chunk}

\begin{thm}\label{fullonpis}
Let $R\to S$ be a map of rings such that $S$ is a finitely presented right $R$-module of finite flat dimension and $\msf{GFlat}(R)$ and $\msf{GFlat}(S)$ are definable. If $S\otimes_{R}-\colon\stgf{R}\to\stgf{S}$ is full and both (C1) and (C2) hold, then $S\otimes_{R}-\colon\msf{GFlat}(R)\to\msf{GFlat}(S)$ is full on pure injective objects.
\end{thm}

\begin{proof}
Let $M$ and $N$ be pure injective Gorenstein flat $R$-modules. Suppose that $f\in\Hom_{S}(S\otimes_{R}M,S\otimes_{R}N)$, and consider $\bar{f}=\pi(f)\colon S\otimes_{R}M\to S\otimes_{R}N$ in $\stgf{S}$. By assumption, there is then a map $\bar{g}\colon M\to N$ in $\stgf{R}$ such that $S\otimes_{R}\bar{g}=\bar{f}$. Let $g\in\Hom_{R}(M,N)$ be a preimage of $\bar{g}$ under $\pi$. The commutating diagram
\[
\begin{tikzcd}
\gfc{R}\arrow[r, "S\otimes_{R}-"] \arrow[d, "\pi"] & \gfc{S} \arrow[d, "\pi"] \\
\stgf{R}\arrow[r, "S\otimes_{R}-"] & \stgf{S}
\end{tikzcd}
\]
ensures that $\bar{f}= \pi(S\otimes_{R}g)$, and thus $S\otimes_{R}g=f+\varphi$ for some $\varphi\colon S\otimes_{R}M\to S\otimes_{R}N$ in $\Mod{S}$ which factors through an object of $\msf{FlatCot}(S)$.

Now, applying \cref{factorsthroughimage}, we may assume that $\phi$ factors through $S\otimes_{R}X$ for some $X\in\msf{FlatCot}(R)$, and thus there are maps $\alpha\colon S\otimes_{R}M\to S\otimes_{R}X$ and $\beta\colon S\otimes_{R}X\to S\otimes_{R}N$ such that $\beta\circ\alpha = \varphi$. The assumptions (C1) and (C2) tell us that $\alpha$ and $\beta$ have preimages under $S\otimes_{R}-$ in $\Hom_{R}(M,X)$ and $\Hom_{R}(X,N)$ respectively, call them $\tilde{\alpha}$ and $\tilde{\beta}$.

Let $\gamma=\tilde{\beta}\circ\tilde{\alpha}\in\Hom_{R}(M,N)$. By construction, $S\otimes_{R}\gamma=\varphi$. It therefore follows that $S\otimes_{R}(g-\varphi)=f$, hence $S\otimes_{R}-$ is full on pure injectives.
\end{proof}

As mentioned in \cref{maponziegler}, the condition of fullness on pure injectives is precisely what is needed to transport indecomposable pure injectives. In particular, one obtains maps on the Ziegler spectrum. Let us make this more precise.

\begin{cor}\label{inducedmaps}
With the assumptions of \cref{fullonpis} holding, yielding that $S\otimes_{R}-\colon\msf{GFlat}(R)\to\msf{GFlat}(S)$ is full on pure injective objects, we have the following:
\begin{enumerate}
\item $S\otimes_{R}-$ preserves pure injective envelopes: $S\otimes_{R}PE(X)=PE(S\otimes_{R}X)$ for all $X\in\msf{GFlat}(R)$;
\item if $X\in\msf{GFlat}(R)$ is indecomposable pure injective, then $S\otimes_{R}X$ is either zero, or is an indecomposable pure injective Gorenstein flat $S$-module;
\item let $\msf{Ker}=\{X\in\msf{GFlat}(R)\cap\msf{pinj}(R):S\otimes_{R}X=0\}$. Then $S\otimes_{R}-$ induces a homeomorphism $\msf{Zg}(\msf{GFlat}(R))\setminus\msf{Ker}\to \msf{Zg}(\msf{pure}(\mc{I}))$, where $\mc{I}$ is the essential image of $S\otimes_{R}-$.
\end{enumerate}
\end{cor}

These corollaries are applications of more general results, which can be found, together with their proofs, at \cite[\S13]{dac}.

Having shown what properties the conditions (C1) and (C2) of \cref{fullconditions} endow the functor $S\otimes_{R}-$, we now illustrate a common situation in which both these conditions hold.

\begin{prop}\label{conditionsholdforcommutative}
Let $R\to S$ be a surjective homomorphism of right coherent rings, such that $S\in\mod{R^{\op}}$ has finite flat dimension. Then both (C1) and (C2) are satisfied.
\end{prop}

\begin{proof}
Viewing $R$ and $S$ as $R-R$-bimodules, there is a short exact sequence of bimodules given by 
\begin{equation}\label{syzygy}
0\to\Omega_{R}^{1}(S)\to R\to S\to 0, 
\end{equation}
where $\Omega_{R}^{1}(S)$ is the first syzygy of $S$. 

Let us first show that (C1) holds. Suppose that $X\in\msf{FlatCot}(R)$ and that $M\in\msf{GFlat}(R)$ is pure injective. Since $X$ is flat, tensoring the above short exact sequence yields an exact sequence of left $R$-modules
\begin{equation}\label{syzses}
0\to \Omega_{R}^{1}(S)\otimes_{R}X\to X\to S\otimes_{R}X\to 0.
\end{equation}
Since $M$ is Gorenstein flat, and $S$ is of finite flat dimension, we may again apply \cite[Lemma 2.3]{cfh} to obtain an isomorphism $S\otimes_{R}M\simeq S\otimes_{R}^{\msf{L}}M$, and we also have $S\otimes_{R}X\simeq S\otimes_{R}^{\msf{L}}X$ as $X$ is flat. There are then isomorphisms
\begin{align*}
\t{Ext}_{S}^{i}(S\otimes_{R}M,S\otimes_{R}X) & = \mrm{H}_{-i}\msf{R}\Hom_{S}(S\otimes_{R}^{\msf{L}}M,S\otimes_{R}^{\msf{L}}X)\\
&\simeq \mrm{H}_{-i}\msf{R}\Hom_{R}(M,\msf{R}\Hom_{S}(S,S\otimes_{R}^{\msf{L}}X))\\
&\simeq\mrm{H}_{-i}\msf{R}\Hom_{R}(M,S\otimes_{R}^{\msf{L}}X)\\
&=\t{Ext}_{R}^{i}(M,S\otimes_{R}X), 
\end{align*}
where the first and last equalities are by the definition of Ext, the second isomorphism is the derived adjunction isomorphism \cite[Theorem 7.5.32]{dcmca}, and the third is the derived counitor \cite[Proposition 7.5.8]{dcmca}.  In particular, by applying $\Hom_{R}(M,-)$ to the short exact sequence in \cref{syzses}, we see that $\Hom_{R}(M,X)\to\Hom_{R}(M,S\otimes_{R}X)$ is surjective provided $\t{Ext}_{R}^{1}(M,\Omega_{R}^{1}(S)\otimes_{R}X)=0$, so this is what we show.

To this end, observe that $S\otimes_{R}X$ is an object of $\msf{GCot}(R)$. Indeed, as $S\otimes_{R}X\in\msf{FlatCot}(S)\subseteq\msf{GCot}(S)$, we have that $\t{Ext}_{S}^{i}(S\otimes_{R}N,S\otimes_{R}X)=0$ for any Gorenstein flat $R$-module $N$, as $S\otimes_{R}-\colon\msf{GFlat}(R)\to\msf{GFlat}(S)$. Thus, by the above isomorphisms, we see that $\t{Ext}_{R}^{i}(N,S\otimes_{R}X)=0$ for any Gorenstein flat $R$-module, as desired. In particular, as $X\in\msf{FlatCot}(R)\subseteq\msf{GCot}(R)$, it follows that the short exact sequence \cref{syzses} is a $\msf{GCot}$-resolution of $\Omega_{R}^{1}(S)\otimes_{R}X$.

Thus, the Gorenstein cotorsion dimension of $\Omega_{R}^{1}(S)\otimes_{R}X$ is at most one. If it is itself Gorenstein cotorsion, then it is immediately clear that $\t{Ext}_{R}^{1}(M,\Omega_{R}(S)\otimes X)=0$. So assume it is not Gorenstein cotorsion. Then there is a minimal Gorenstein cotorsion resolution
\begin{equation}\label{gcotres}
0\to\Omega_{R}^{1}(S)\otimes_{R}X\to C^{0}\to C^{1}\to 0,
\end{equation}
where $\Omega_{R}^{1}(S)\otimes_{R}X\to C^{0}$ is a special $\msf{GCot}$-envelope, which exists as $(\msf{GFlat}(R),\msf{GCot}(R))$ is a complete cotorsion pair, as detailed in \cref{GFlatcotorsionpair}. In particular, it follows that $C^{1}\in\msf{GFlat}(R)$ by Salce's lemma in \cref{cotorsionpair}. Consequently, $C^{1}\in\msf{GFlat}(R)\cap\msf{GCot}(R)=\msf{FlatCot}(R)$, and thus $C^{1}$ is flat.

Now, as $\Omega_{R}^{1}(S)$ is a finitely presented $R-R$-bimodule by assumption, the functor $\Omega_{R}^{1}(S)\otimes_{R}-\colon\Mod{R}\to\Mod{R}$ is a definable functor, and thus preserves pure injective objects. Yet, by assumption $X$ was pure injective, hence $\Omega_{R}^{1}(S)\otimes_{R}X\in\msf{Pinj}(R)$, and is therefore cotorsion by \cref{Flatcotorsionpair}. In particular, we see that
\[
\t{Ext}_{R}^{1}(C^{1},\Omega_{R}^{1}(S)\otimes_{R}X)=0,
\]
as $(\msf{Flat}(R),\msf{Cot}(R))$ is a cotorsion pair, and therefore \cref{gcotres} splits. Thus, $\Omega_{R}^{1}(S)\otimes_{R}X$ is a summand of $C^{0}$, which is Gorenstein cotorsion, and $\Omega_{R}^{1}(S)\otimes_{R}X$ is itself Gorenstein cotorsion, meaning $\t{Ext}_{R}^{1}(M,\Omega_{R}^{1}(S)\otimes_{R}X)=0$. This finishes the proof that (C1) holds.

The proof for the second condition is very similar. Tensoring \cref{syzygy} with $M$, and noting that $\t{Tor}_{i}^{R}(S,M)=0$ for all $i>0$, the same line of reasoning as above shows that it suffices to prove that $\t{Ext}_{R}^{1}(X,\Omega_{R}^{1}(S)\otimes_{R}M)=0$. Yet this happens immediately from the facts that $X\in\Flat{R}$, and $\Omega_{R}^{1}(S)\otimes_{R}M\in\msf{Pinj}(R)$.
\end{proof}

\begin{chunk}\label{kcdefinition}
With the necessary assumptions such that $R\to S$ induces a definable functor $\msf{GFlat}(R)\to\msf{GFlat}(R)$ of definable subcategories, the results of \cref{inducedmaps} illustrate that understanding the kernel of this map is quite a significant part of determining the transport of information in the Ziegler spectrum.

The kernel of a definable functor is always a definable subcategory (being the preimage of the zero definable subcategory), which is either zero or contains a non-zero cotorsion module; in either case, we only need again consider
\[
\msf{K}=\{X\in\gfc{R}:S\otimes_{R}X=0\}.
\]

By \cref{triangulatedfunctor}, we also obtain a product preserving triangulated functor $\stgf{R}\to\stgf{S}$, whose kernel is thus a colocalising subcategory of $\stgf{R}$; in fact, when the categories involved are compactly generated and $S\otimes_{R}-$ also preserves coproducts, this kernel is also a definable subcategory of $\stgf{R}$. We shall let
\[
\msf{C}=\{X\in\stgf{R}:S\otimes_{R}X=0\}
\]
denote this subcategory.

We now aim to relate $\msf{K}$ with $\msf{C}$. To do this, let us observe that $\msf{K}$ has the structure of an exact category, as $S\otimes_{R}-\colon\msf{GFlat}(R)\to\msf{GFlat}(S)$ is an exact functor; in particular $\msf{K}\cap\Flat{R}$ is also an exact subcategory of $\msf{K}$.
\end{chunk}

\begin{prop}\label{kernelfrobenius}
Let $R\to S$ be a morphism of right coherent rings such that $S\in\mod{R^{\op}}$ is of finite flat dimension. If, for every $X\in\msf{K}$ there is a deflation $F\to X$ and an inflation $X\to G$, with $F,G\in\msf{K}\cap\Flat{R}$, then $\msf{K}$ is a Frobenius category with injective objects $\msf{K}\cap\msf{Flat}(R)$, whose stable category $\ul{\msf{K}}$ is a colocalising subcategory of $\stgf{R}$ contained in $\msf{C}$.
\end{prop}

\begin{proof}
Let us first show that an object of $\msf{K}$ is projective and injective if and only if it lies in $\msf{K}\cap\Flat{R}$. Clearly, since $\msf{K}\cap\Flat{R}\subseteq\msf{FlatCot}(R)$, any object in $\msf{K}\cap\Flat{R}$ is projective and injective. For the converse, suppose that $P\in\msf{K}$ is a projective object. By assumption there is an exact sequence $0\to P'\to F\to P\to 0$ with $F\in\msf{K}\cap\Flat{R}$, and as $P$ is projective this sequence splits; thus $P$ is flat. A dual argument using an exact sequence $0\to I\to G\to I'\to 0$, where $G$ is flat, yields that an injective object is also flat. In particular, $\msf{K}$ also has enough projectives and injectives.

Let us now show that $\ul{\msf{K}}$ is colocalising in $\stgf{R}$.
Using the diagram in \cref{diagram}, we see that any object in $\ul{\msf{K}}$ is also in $\msf{C}$; moreover, as any object of $\msf{K}\cap\Flat{R}$ is in $\msf{FlatCot}(R)$, there is an inclusion of homomorphisms $\underline{\Hom}_{\msf{K}}(X,Y)\subseteq \Hom_{\msf{C}}(X,Y)$ for any pair of objects $X,Y\in\msf{K}$. For the reverse inclusion, suppose that $X\to Y$ factors through a flat cotorsion module $F$. By assumption, there is an inflation $0\to X\to G$ with $G\in\msf{K}\cap\Flat{R}$, and as $F$ is injective in $\gfc{R}$ the map $X\to F$ factors through $G$, hence $X\to Y$ factors through an object in $\msf{K}\cap\Flat{R}$. This shows that $\ul{\msf{K}}$ is a full subcategory of $\stgf{R}$.

To finish the proof, we must show that $\ul{\msf{K}}$ is closed under products and triangles in $\stgf{R}$. The former is straightforward as $\msf{K}$ is closed under products. For the latter, a triangle $X\to Y\to Z\to \Sigma X$ in $\stgf{R}$ arises from a short exact sequence $0\to X\to Y\to Z\to 0$ in $\gfc{R}$ as detailed in \cite[\S 3.3]{krbook}. As $S\otimes_{R}-\colon\gfc{R}\to\gfc{S}$ is an exact functor, we see that if any two of $X,Y$ and $Z$ lie in $\msf{K}$, then the other does as well, which shows that $\ul{\msf{K}}$ is triangulated with respect to the triangulation on $\stgf{R}$.
\end{proof}

A natural question that arises from \cref{kernelfrobenius} is how different are the objects of $\ul{\msf{K}}$ and $\msf{C}$? Since we are predominantly interested in purity considerations, it makes sense to try to better establish a relationship between the pure injective objects in both categories. It transpires that when the functor $S\otimes_{R}-$ possesses a particular property, there is no difference between the pure injective objects in $\ul{\msf{K}}$ and $\msf{C}$. Before being more precise, let us introduce the class
\[
\msf{X}=\{X\in\msf{GFlat}(R):S\otimes_{R}X\in\Flat{S}\},
\]
which is a definable subcategory of $\msf{GFlat}(R)$ whenever $\msf{GFlat}(R)$ is definable. It is clear that $\msf{X}$ consists of the Gorenstein flat $R$-modules $X$ such that $S\otimes_{R}X$ vanishes in $\stgf{S}$, and can therefore be thought of as a definable category which measures the `discrepancy' between $\msf{K}$ and $\msf{C}$.

\begin{prop}
Let $R\to S$ be a map of coherent rings such that $S\in\mod{R^{\op}}$ has finite flat dimension, and assume that $\msf{GFlat}(R)$ and $\msf{GFlat}(S)$ are definable. If the restriction $S\otimes_{R}-\colon\msf{X}\to\Flat{S}$ is full on pure injectives, then there are equalities
\[
\msf{pinj}(\msf{C})= \msf{pinj}(\msf{X}\setminus \Flat{R}) = \msf{pinj}(\msf{K}\setminus\Flat{R}).
\]
In particular, the indecomposable pure injective objects in $\ul{\msf{K}}$ and $\msf{C}$ are the same.
\end{prop}

\begin{proof}
Let us first show that an indecomposable pure injective object is in $\msf{C}$ if and only if it is in $\msf{X}$. If $X\in\msf{X}\setminus\Flat{R}$ is an indecomposable pure injective object, then, as the image of $S\otimes_{R}X$ in $\stgf{S}$ is zero, it follows from \cref{diagram} that either $X$ is zero in $\stgf{R}$, which is not possible by assumption, or $X\in\msf{C}$ when viewed as an object of $\stgf{R}$. On the other hand, if $Z\in\msf{C}$ is an indecomposable pure injective, then, by \cref{piinst}, $Z$ is also indecomposable pure injective in $\gfc{S}$. By the commutativity of \cref{diagram}, it follows that $S\otimes_{R}Z$ is a flat $S$-module, hence $Z\in\msf{X}$.

It is clear that any indecomposable pure injective object in $\msf{K}$ that is not flat is in $\msf{X}\setminus\Flat{R}$, so the proof is complete once the converse inclusion is shown. By the assumption that $\msf{X}\to\Flat{S}$ is full on pure injectives, there is a homeomorphism $\msf{Zg}(\msf{X})\setminus\{X\in\msf{Zg}(\msf{X}):S\otimes_{R}X=0\}\simeq \msf{Zg}(\Flat{S}),$
where the codomain is $\msf{Zg}(\Flat{S})$ by \cref{psoofimage}. However, $\{X\in\msf{Zg}(\msf{X}):S\otimes_{R}X=0\}$ is nothing other than the set of indecomposable pure injective objects in $\msf{K}$. In particular, there is a homeomorphism
\begin{equation}\label{homeo}
\msf{Zg}(\Flat{S})\simeq \msf{Zg}(\msf{X})\setminus \msf{Zg}(\msf{K}).
\end{equation}
However, we also have that the restriction $S\otimes_{R}-\colon\Flat{R}\to\Flat{S}$ is full on pure injectives, hence induces a homeomorphism $\msf{Zg}(\Flat{S})\simeq \msf{Zg}(\Flat{R})\setminus\msf{Zg}(\Flat{R}\cap\msf{K})$, which is just a restriction of \cref{homeo}. In particular, any indecomposable pure injective module in $\msf{X}\setminus\msf{K}$ is a flat module; consequently $\msf{pinj}(\msf{X}\setminus\Flat{R})=\msf{pinj}(\msf{K})$, which is what was required.
\end{proof}

\section{Phenomena over commutative rings}\label{commrings}
In this section we focus on two particular phenomena that happen over commutative Noetherian rings, and we use tools which are specific to the commutative setting. The first of the phenomena will enable us to give an example where the conditions of \cref{kernelfrobenius} are satisfied, while in the second we exploit the relationship between Gorenstein flat modules and maximal Cohen-Macaulay modules over hypersurfaces to develop an analogy for Kn\"{o}rrer periodicity in the setting of Gorenstein flat cotorsion modules. 

\subsection{Regular sequences and closure under flat covers and cotorsion envelopes}

Let $R$ be a commutative Noetherian ring, and suppose that $\bm{x}$ is an $R$-sequence. The canonical ring map $R\to R/\bm{x}$ satisfies the main assumptions we have considered so far; in this case, as $R$ is Noetherian, all that matters is that $R/\bm{x}$ has finite flat dimension over $R$, which was discussed in \cref{regseq}. 

Thus, let $I=(\bm{x})$ be an ideal of $R$, where $\bm{x}=x_{1},\cdots,x_{n}$ is a regular sequence, and consider the functor $R/I\otimes_{R}-\colon\msf{GFlat}(R)\to\msf{GFlat}(R/I)$ viewed along the quotient $R\to R/I$. Recall from \cref{kcdefinition} the associated categories $\msf{K}$ and $\msf{C}$. The goal of this section is to prove the following theorem.

\begin{thm}\label{comkerfrob}
Let $R$ be a commutative Noetherian ring and let $I$ be an ideal generated by an $R$-sequence. Then $\msf{K}=\{X\in\gfc{R}:S\otimes_{R}X=0\}$ is closed under flat covers and Gorenstein cotorsion (pre)envelopes. In particular, the conditions of \cref{kernelfrobenius} are satisfied, and $\msf{K}$ is a Frobenius category.
\end{thm}

Before proving this result, let us recall some further background material from commutative algebra.

\begin{chunk}
Suppose that $\mf{a}$ is an ideal of a commutative Noetherian ring $R$. We can associate two functors to $\mf{a}$, the $\mf{a}$-\ti{torsion} and $\mf{a}$-\ti{adic completion} functors; which are, respectively, defined to be
\[
\Gamma_{\mf{a}}(-):=\rlim_{i}\Hom_{R}(R/\mf{a}^{i},-)
\]
and
\[
\Lambda^{\mf{a}}(-):=\llim_{i}R/\mf{a}^{i}\otimes_{R}-.
\]
These induce derived functors
\[
\msf{R}\Gamma_{\mf{a}}(-)\colon\D(R)\to\D(R)
\]
and 
\[
\msf{L}\Lambda^{\mf{a}}(-)\colon\D(R)\to \D(R),
\]
which are called derived torsion and derived completion respectively. The relationship between derived torsion and derived completion is given by Greenlees-May duality, which states that
\begin{equation}\label{greenleesmay}
\msf{R}\Hom_{R}(\msf{R}\Gamma_{\mf{a}}(X),Y)\simeq \msf{R}\Hom_{R}(X,\msf{L}\Lambda^{\mf{a}}(Y))
\end{equation}
for any ideal $\mf{a}$ and any complexes $X,Y\in\D(R)$.

If $M$ is an $R$-module, the $i$-th \ti{local cohomology of $M$ with support in $\mf{a}$} is defined to be
\[
H_{\mf{a}}^{i}(M)=\mrm{H}_{-i}\msf{R}\Gamma_{\mf{a}}(M).
\]
A proof of the above adjunction (in its most general form), as well as some historical information about it, can be found in \cite{psy}.
\end{chunk}

\begin{chunk}\label{commflatcotorsion}
Over commutative Noetherian rings, modules in $\msf{FlatCot}(R)$ have a particularly pleasant structure. Recall, from for example \cite[Proposition 3.1]{matlis}, over such a ring $R$, given a prime $\mf{p}$, the injective hull $E(R/\mf{p})$ of $R/\mf{p}$ is indecomposable, and every indecomposable arises this way. Since $R$ is Noetherian, any coproduct of injective modules is injective, and thus the module
\[
\Hom_{R}(E(R/\mf{p}),E(R/\mf{p})^{(X_{\mf{p}})})
\]
is flat and cotorsion for any indexing set $X_{\mf{p}}$. Such a module can, by Matlis duality, be seen as the completion of a free $R_{\mf{p}}$-module of rank $\t{card}(X_{\mf{p}})$, by \cite[Thereom 3.4.1]{rha}. In fact, as shown in \cite[Theorem 5.3.28]{rha}, if $F\in\msf{FlatCot}(R)$, then
\begin{equation}\label{fcdecomp}
F\simeq\prod_{\mf{p}\in\mrm{Spec}(R)}\Hom_{R}(E(R/\mf{p}),E(R/\mf{p})^{(X_{\mf{p}})})
\end{equation}
for some sets $X_{\mf{p}}$. 

If $M\in\msf{Cot}(R)$ then $F(M)$, the flat cover of $M$, lies in $\msf{FlatCot}(R)$, by \cref{Flatcotorsionpair}. In particular, it has a decomposition as in \cref{fcdecomp}; the cardinality of the set $X_{\mf{p}}$ in the decomposition of $F(M)$ is denoted by $\pi_{0}(\mf{p},M)$, which was shown in \cite[Theorem 5.2.2]{xu} to be equal to
\[
\pi_{0}(\mf{p},M)=\t{dim}_{k(\mf{p})}(k(\mf{p})\otimes_{R_{\mf{p}}}\Hom_{R}(R_{\mf{p}},M)),
\]
where $k(\mf{p})$ is the residue field of $R_{\mf{p}}$. For brevity, we write $T_{\mf{p}}=\Hom_{R}(E(R/\mf{p}),E(R/\mf{p})^{(X_{\mf{p}})})$ for this completion of a free $R_{\mf{p}}$-module.
\end{chunk}

Let us now proceed with the proof of \cref{comkerfrob}.

\begin{proof}[Proof of \cref{comkerfrob}]
Let us first show that $\msf{K}$ is closed under Gorenstein cotorsion pre-envelopes. Suppose that $M\in\msf{K}$, and consider the short exact sequence
\[
0\to M\to X\to N\to 0,
\]
where $M\to X$ is a Gorenstein cotorsion pre-envelope and $N$ is a Gorenstein flat $R$-module. This short exact sequence exists by \cref{GFlatcotorsionpair}. Since $M$ is itself Gorenstein flat, and $\msf{GFlat}(R)$ is closed under extensions, $X$ is in $\msf{GFlat}(R)\cap\msf{GCot}(R)=\msf{FlatCot}(R)$; moreover $N$ is a cotorsion module as both $M$ and $X$ are, since any Gorenstein cotorsion module is cotorsion.

As $X\in\msf{FlatCot}(R)$, it has a decomposition
\[
X\simeq \prod_{\mf{p}\in\mrm{Spec}(R)}T_{\mf{p}};
\]
thus $\Hom_{R}(M,X)\simeq\prod_{\mf{p}\in\t{Spec}(R)}\Hom_{R}(M,T_{\mf{p}})$. Now, $M\in\gfc{R}$, and $T_{\mf{p}}\in\msf{FlatCot}(R)\subseteq\msf{GCot}(R)$, hence $\t{Ext}_{R}^{i}(M,T_{\mf{p}})=0$ for all $i>0$, giving an isomorphism $\Hom_{R}(M,T_{\mf{p}})\simeq\msf{R}\Hom_{R}(M,T_{\mf{p}})$ in $\D(R)$. Now, as $T_{\mf{p}}$ is flat and is the completion of a free $R_{\mf{p}}$-module (in the $\mf{p}$-adic topology), we have $T_{\mf{p}}\simeq \msf{L}\Lambda^{\mf{p}}(\oplus_{X_{\mf{p}}}R_{\mf{p}})$ by \cite[Corollary 13.1.16]{dcmca}. Combining these facts, we get isomorphisms
\begin{align}
\msf{R}\Hom_{R}(M,T_{\mf{p}})&\simeq \msf{R}\Hom_{R}(M,\msf{L}\Lambda^{\mf{p}}(\oplus_{X_{\mf{p}}}R_{\mf{p}})) \nonumber\\
&\simeq \msf{R}\Hom_{R}(\msf{R}\Gamma_{\mf{p}}M,\oplus_{X_{\mf{p}}}R_{\mf{p}}) \mbox{ by \cref{greenleesmay}} \label{gmisos}.
\end{align}
Now, the sequence $\bm{x}$ is either contained in $\mf{p}$ or it is not. If it is, then by assumption of $M$ being an object of $\msf{K}$, it follows that $R/(\bm{x})\otimes_{R}^{\msf{L}}M=0$, and consequently $\msf{R}\Hom_{R}(R/(\bm{x}),M)=0$ as well. In particular, it follows that $\msf{R}\Hom_{R}(R/\mf{p},M)=0$ by \cite[Proposition 5.3.11]{strooker}. We may therefore apply \cite[Proposition 5.3.15]{strooker} to see that $H_{\mf{p}}^{i}(M)=0$ for all $i\geq 0$, and consequently $\msf{R}\Gamma_{\mf{p}}(M)\simeq 0$.

Using this fact, and \cref{gmisos}, it follows that $\msf{R}\Hom_{R}(M,T_{\mf{p}})$ is acyclic, and thus $\Hom_{R}(M,T_{\mf{p}})=0$ whenever $\bm{x}\subseteq\mf{p}$. Consequently we have that 
\[
X\simeq\prod_{\bm{x}\not\subseteq\mf{p}}T_{\mf{p}}.
\]

So, suppose that $\bm{x}\not\subseteq\mf{p}$, which, by \cite[Lemma 3.2.9]{strooker}, tells us that $\msf{R}\Hom_{R}(R/\bm{x},E(R/\mf{p}))\simeq 0$, and thus $R/\bm{x}\otimes_{R}^{\msf{L}}E(R/\mf{p})\simeq 0$, by \cite[Theorem 14.4.12]{dcmca}. Consequently, 
\[
\msf{R}\Hom_{R}(R/(\bm{x}),T_{\mf{p}})\simeq \msf{R}\Hom_{R}(R/(\bm{x})\otimes_{R}^{\msf{L}} E(R/\mf{p}),E(R/\mf{p})^{(X_{\mf{p}})})\simeq 0,
\]
and therefore, by another application of \cite[Theorem 14.4.12]{dcmca}, we have that $R/(\bm{x})\otimes_{R}^{\msf{L}}T_{\mf{p}}=R/(\bm{x})\otimes_{R} T_{\mf{p}}=0$. Therefore
\[
R/(\bm{x})\otimes X\simeq \prod_{\bm{x}\not\subseteq\mf{p}}R/(\bm{x})\otimes T_{\mf{p}}=0,
\]
so $X\in\msf{K}$ as desired.

We now show that $\msf{K}$ is closed under flat covers. Suppose that $F\to M$ is the flat cover of $M$, so $F$ is cotorsion. Again, express $F\simeq\prod T_{\mf{p}}$. As stated in \cref{commflatcotorsion}, $T_{\mf{p}}$ is determined by
\[
\pi_{0}(\mf{p},M)=\t{dim}_{k(\mf{p})}(k(\mf{p})\otimes_{R_{\mf{p}}}\Hom_{R}(R_{\mf{p}},M)). 
\]
We showed above that $R/\bm{x}\otimes_{R}T_{\mf{p}}\neq 0$ if and only if $\bm{x}\subseteq\mf{p}$, hence it suffices to show $\pi_{0}(\mf{p},M)=0$ for all primes $\mf{p}$ containing $\bm{x}$.

As $M$ is cotorsion and $R_{\mf{p}}$ is flat over $R$, we have that $\Hom_{R}(R_{\mf{p}},M)\simeq\msf{R}\Hom_{R}(R_{\mf{p}},M)$, hence if $k(\mf{p})\otimes_{R}^{\msf{L}}\msf{R}\Hom_{R}(R_{\mf{p}},M)$ is acyclic, it follows that $\pi_{0}(\mf{p},M)=0$, and this is what we show.

By \cite[Proposition 4.4]{sww}, the complex $k(\mf{p})\otimes_{R}^{\msf{L}}\msf{R}\Hom_{R}(R_{\mf{p}},M)$ is acyclic if and only if $\msf{R}\Hom_{R}(k(\mf{p}),M)$ is acyclic. As $k(\mf{p})\simeq R_{\mf{p}}\otimes_{R}^{\msf{L}}R/\mf{p}$, we have
\[
\msf{R}\Hom_{R}(k(\mf{p}),M)\simeq \msf{R}\Hom_{R}(R_{\mf{p}}\otimes_{R}^{\msf{L}}R/\mf{p},M)\simeq \msf{R}\Hom_{R}(R_{\mf{p}},\msf{R}\Hom_{R}(R/\mf{p},M)).
\]
Since $M\in\msf{K}$, we have $R/\bm{x}\otimes_{R}M=0$, and thus $\t{Tor}_{i}^{R}(R/\bm{x},M)=0$ for all $i\geq 0$, as $R/\bm{x}$ has finite flat dimension over $R$; hence $R/\bm{x}\otimes_{R}^{\msf{L}}M\simeq 0$. But as $\bm{x}\subset\mf{p}$, it follows that $R/\mf{p}\otimes_{R}^{\msf{L}}M\simeq 0$ which is equivalent to  $\msf{R}\Hom_{R}(R/\mf{p},M) \simeq 0$ by \cite[Theorem 14.4.12]{dcmca}. Consequently $\msf{R}\Hom_{R}(k(\mf{p}),M)=0$, and thus $k(\mf{p})\otimes_{R}^{\msf{L}}\msf{R}\Hom_{R}(R_{\mf{p}},M)$ is acyclic as desired. 
\end{proof}

\begin{rem}
Note that we only needed $(\bm{x})$ to be regular in the proof that $\msf{K}$ is closed under flat covers. The above proof actually shows that if $I$ is any ideal, then $\{X\in\gfc{R}:X=IX\}$ is closed under Gorenstein cotorsion pre-envelopes. 
\end{rem}

\subsection{The double-double branched cover and Kn\"{o}rrer periodicity for Gorenstein flat modules}\label{section:knorrer}
Let us now consider the particular case of Gorenstein flat modules over hypersurfaces, in which we can establish an extension of Kn\"{o}rrer periodicity. Let us first give some specific background information.

\begin{chunk}
Let $k$ be an algebraically closed field of characteristic other than two. Let $S=k[[x_{1},\cdots,x_{n+1}]]$ denote the power series ring in $n+1$ variables. This is a complete local ring whose maximal ideal shall be denoted $\mf{n}$; moreover, it is a domain, hence any non-zero element of $S$ is a regular element. Given any non-zero $f\in\mf{n}^{2}$, we let $R=S/(f)$, which is a Gorenstein local ring of Krull dimension $n$. We shall let $\mf{m}$ denote the maximal ideal of $R$. 

Over Gorenstein rings, the category of Gorenstein flat modules is a finitely accessible definable subcategory; its finitely presented modules are precisely the finitely presented Gorenstein projective modules. Over Gorenstein local rings, the category of finitely presented Gorenstein projective modules coincides with the \ti{maximal Cohen-Macaulay} modules. A finitely presented $R$-module $M$ is maximal Cohen-Macaulay if
\[
H_{\mf{m}}^{i}(M)\neq 0
\]
if and only if $i=n=\t{dim}\,R$. We let $\msf{CM}(R)$ denote the full subcategory of $\mod{R}$ consisting of maximal Cohen-Macaulay modules, and (as stated) $\rlim\msf{CM}(R)=\msf{GFlat}(R)$. In particular $\msf{CM}(R)$ is a Frobenius category.
\end{chunk}

\begin{chunk}
The \ti{double branched cover} of $R$ is the ring 
\[
R^{\sharp}:=S[[z]]/(f+z^{2}),
\]
which is clearly also a hypersurface. The element $z\in R^{\sharp}$ is regular, and the quotient ring $R^{\sharp}/(z)$ is isomorphic to $R$. Moreover, extending scalars along $R^{\sharp}\to R$ gives a functor
\[
\msf{CM}(R^{\sharp})\to\msf{CM}(R),
\]
which is a functor of Frobenius categories, and the induced functor on the corresponding singularity categories is triangulated. 
\end{chunk}

\begin{chunk}
A \ti{matrix factorisation} of $f$ is a pair of square matrices $(\phi,\psi)$ of the same finite size over $S$ such that $\phi\psi=\psi\phi=f\t{Id}$. The category of matrix factorisations of $f$ is denoted $\msf{MF}(S,f)$. There are two distinguished matrix factorisations, namely $(1,f)$ and $(f,1)$. Given a matrix factorisation $(\phi,\psi)$, the module $\msf{Cok}(\phi,\psi)=\msf{Cok}(\phi\colon S^{(n)}\to S^{(n)})$ is a maximal Cohen-Macaulay $R$-module.

A classic result of Eisenbud states that $\msf{CM}(R)$ is equivalent to the category $\msf{MF}(S,f)/\{(1,f)\}$, see \cite[Theorem 8.7]{lw}. Furthermore, there is actually a triangulated equivalence $\ul{\msf{CM}}(R)\simeq \msf{MF}(S,f)/\{(1,f),(f,1)\}$. 

In \cite{knorrer}, Kn\"{o}rrer used matrix factorisations, among other things, to relate the categories $\ul{\msf{CM}}(R)$ and $\ul{\msf{CM}}(R^{\sharp\sharp})$. He constructed a functor $(-)^{\dagger}\colon\msf{CM}(R)\to\msf{CM}(R^{\sharp\sharp})$ of Frobenius categories which induced a triangulated functor 
\begin{equation}\label{knorequivalence}
(-)^{\dagger}\colon\ul{\msf{CM}}(R)\to\ul{\msf{CM}}(R^{\sharp\sharp}),
\end{equation}
which, by \cite[Theorem 3.1]{knorrer}, is an equivalence of categories.
\end{chunk}

Let us first prove that Kn\"{o}rrer's equivalence extends to an equivalence of categories between the stable categories of Gorenstein flat cotorsion modules.

\begin{thm}\label{bigknorrer}
Let $R=S/(f)$ be a hypersurface singularity, and let $R^{\sharp\sharp}=S[[z_{1},z_{2}]]/(f+z_{1}^{2}+z_{2}^{2})$. Then Kn\"{o}rrer's equivalence of \cref{knorequivalence} extends to a triangulated equivalence \[\stgf{R}\to\stgf{R^{\sharp\sharp}}.\]
\end{thm}

\begin{proof}
As shown in \cite[Definiton 2.1]{dyckerhoff}, there is a dg-category of matrix factorisations of $f$ over $S$, which we denote $\msf{MF}_{\t{dg}}(S,f)$; the category $\msf{MF}(S,f)$ is the category obtained by taking cycles of the morphism complexes of $\msf{MF}_{\t{dg}}(S,f)$, while the homotopy category of $\msf{MF}_{\t{dg}}(S,f)$ is triangle equivalent to $\ul{\msf{CM}}(R)$. 

Combining \cite[\S5.3]{dyckerhoff} with \cite[Proposition 2.20]{brownknorrer}, shows that Kn\"{o}rrer's equivalence of \cref{knorequivalence} arises as the derived functor of a dg functor $\tilde{\dagger}\colon\msf{MF}_{\t{dg}}(S,f)\to\msf{MF}_{\t{dg}}(S[[z_{1},z_{2}]],f+z_{1}^{2}+z_{2}^{2})$. We may, by taking the dg nerve of \cite[1.3.1]{HA}, pass to the $\infty$-categorical setting, and view $\tilde{\dagger}\colon\mathscr{D}_{\t{sg}}(R)\to\mathscr{D}_{\t{sg}}(R^{\sharp\sharp})$ as a functor between the stable $\infty$-categorical enhancements of $\ul{\msf{CM}}(R)$ and $\ul{\msf{CM}}(R^{\sharp\sharp})$ respectively. Since $\dagger\colon\msf{h}\mathscr{D}_{\t{sg}}(R)\to\msf{h}\mathscr{D}_{\t{sg}}(R^{\sharp\sharp})$ is an equivalence, so is $\tilde{\dagger}$. 

By taking ind-completions, which are functorial, we obtain an equivalence $\tilde{\dagger}\colon \msf{Ind}\mathscr{D}_{\t{sg}}(R)\to\msf{Ind}\mathscr{D}_{\t{sg}}(R^{\sharp\sharp})$. 
But $\msf{Ind}\mathscr{D}_{\t{sg}}(R)\simeq \mathscr{K}_{\t{ac}}(\msf{Inj}(R))$, an $\infty$-category enhancement of $\msf{K}_{\t{ac}}(\msf{Inj}(R))$. This is because $\mathscr{K}_{\t{ac}}(\msf{Inj}(R))$ is presentable, by \cite[Corollary 1.4.4.2]{HA}, and presentable infinity categories are the ind-completion of their compact objects by \cite[Theorem 5.5.1.1(6)]{htt}, and the compact objects of $\mathscr{K}_{\t{ac}}(\msf{Inj}(R))$ is $\mathscr{D}_{\t{sg}}(R)$ by \cite{krstable} (see \cite[\S 8.3]{antieau2021uniqueness} as well). 

Consequently, we get that 
\[
\msf{h}(\tilde{\dagger})\colon\msf{K}_{\t{ac}}(\msf{Inj}(R))\to\msf{K}_{\t{ac}}(\msf{Inj}(R^{\sharp\sharp})) 
\]
is an equivalence of categories. The fact that $\stgf{R}\simeq\msf{K}_{\t{ac}}(\msf{Inj}(R))$, and likewise for $R^{\sharp\sharp}$, yields the claim.
\end{proof}

Let us now see how this periodicity can be used to compare the Ziegler spectra of $\msf{GFlat}(R)$ and $\msf{GFlat}(R^{\sharp\sharp})$. Recall from \cref{kcdefinition} that $\msf{K}=\{X\in\gfc{R}:R\otimes_{R^{\sharp\sharp}}X=0\}$.

\begin{prop}\label{hypziegler}
There is a homeomorphism 
\[
\msf{Zg}(\msf{GFlat}(R))\simeq\msf{Zg}(\msf{GFlat}(R^{\sharp\sharp}))\setminus\msf{Zg}(\msf{K}\cap\Flat{R^{\sharp\sharp}}).
\]
The points of $\msf{Zg}(\msf{K}\cap\Flat{R^{\sharp\sharp}})$ are the $T_{\mf{p}}$, the completion of a free $R_{\mf{p}}$-module of rank one, where $\mf{p}\not\in V(z_{1},z_{2})\subseteq\msf{Spec}(R^{\sharp\sharp})$.
\end{prop}

\begin{proof}
By \cref{bigknorrer}, there is clearly an equivalence $\msf{Zg}(\stgf{R})\simeq\msf{Zg}(\stgf{R^{\sharp\sharp}})$. By \cref{zieglerdisjoint}, we have that $\msf{Zg}(\msf{GFlat}(R))\simeq\msf{Zg}(\stgf{R})\coprod\msf{Zg}(\Flat{R})$. Combining these facts, there is a homeomorphism
\begin{equation}\label{zgdqknor}
\msf{Zg}(\msf{GFlat}(R))\setminus\msf{Zg}(\Flat{R})\simeq \msf{Zg}(\msf{GFlat}(R^{\sharp\sharp}))\setminus \msf{Zg}(\Flat{R^{\sharp\sharp}}).
\end{equation}
By restricting, we obtain a functor $R\otimes_{R^{\sharp\sharp}}-\colon\Flat{R^{\sharp\sharp}}\to\Flat{R}$. As illustrated in the proof of \cref{comkerfrob}, we have that, for $T_{\mf{p}}\in\Flat{R}$ pure injective with $\mf{p}$ a prime, $R\otimes_{R^{\sharp\sharp}}T_{\mf{p}}\neq 0$ if and only if $(z_{1},z_{2})\subseteq\mf{p}$. In particular one has $T_{\mf{q}}\in\msf{K}$ if and only if $\mf{q}\not\in V(z_{1},z_{2})$. This provides the criterion for the points in $\msf{Zg}(\msf{K}\cap\Flat{R^{\sharp\sharp}})$. 

Now, the restriction $R\otimes_{R^{\sharp\sharp}}-\colon\Flat{R^{\sharp\sharp}}\to\Flat{R}$ is full on pure injectives by \cref{conditionsholdforcommutative}, and therefore preserves indecomposable pure injective objects. Moreover, by \cref{psoofimage}, every indecomposable pure injective object of $\Flat{R}$ is of the form $R\otimes_{R^{\sharp\sharp}} F$ for some $F\in\Flat{R^{\sharp\sharp}}$, and $F$ is necessarily of the form $T_{\mf{p}}$ for $\mf{p}\in V(z_{1},z_{2})$. In particular, there is a homeomorphism
\[
\msf{Zg}(\Flat{R})\simeq \{T_{\mf{q}}:\mf{q}\in V(z_{1},z_{2})\subseteq\msf{Spec}(R^{\sharp\sharp})\}.
\]
By \cite[Theorem 15.5]{dac}, there is a homeomorphism $\msf{Zg}(\Flat{R^{\sharp\sharp}})\setminus\msf{Zg}(\msf{K}\cap\Flat{R^{\sharp\sharp}})\simeq\msf{Zg}(\Flat{R})$. Combining these facts, and taking the disjoint union of $\msf{Zg}(\Flat{R})$ with both sides of \cref{zgdqknor}, we obtain a homeomorphism
\[
\msf{Zg}(\msf{GFlat}(R))\simeq\msf{Zg}(\msf{GFlat}(R^{\sharp\sharp}))\setminus\msf{Zg}(\msf{K}\cap\Flat{R^{\sharp\sharp}}).
\] 
as desired.
\end{proof}

Let us now investigate further properties of the functor $R\otimes_{R^{\sharp\sharp}}-$ in terms of adjoints it admits. The first thing we need is the following lemma, which follows quickly from Kn\"{o}rrer's work.

\begin{lem}\label{functorsoncompacts}
The functor $R\otimes_{R^{\sharp\sharp}}-\colon\ul{\msf{CM}}(R^{\sharp\sharp})\to\ul{\msf{CM}}(R)$ has a left and right adjoint, which is given by $f:=(-\oplus\Sigma-)^{\dagger}\colon \ul{\msf{CM}}(R)\to\ul{\msf{CM}}(R^{\sharp\sharp})$.
\end{lem}

\begin{proof}
Suppose $M\in\ul{\msf{CM}}(R^{\sharp\sharp})$ and $N\in\ul{\msf{CM}}(R)$; then there are isomorphisms
\begin{align*}
\Hom_{\ul{\msf{CM}}(R)}(R\otimes_{R^{\sharp\sharp}}M,N) &\simeq \Hom_{\ul{\msf{CM}}(R)}(R\otimes_{R^{\sharp\sharp}}A^{\dagger},N) \t{ for some $A\in\ul{\msf{CM}}(R)$ by \cref{knorequivalence}}\\
&\simeq \Hom_{\ul{\msf{CM}}(R)}(A\oplus\Sigma^{-1}A,N), \t{ by \cite[Lemma 8.32]{lw}},\\
&\simeq \Hom_{\ul{\msf{CM}}(R)}(A,N\oplus \Sigma N)\\
&\simeq \Hom_{\ul{\msf{CM}}(R)}(A, \Sigma(N\oplus \Sigma^{-1}N))\\
&\simeq \Hom_{\ul{\msf{CM}}(R^{\sharp\sharp})}(A^{\dagger}, N^{\dagger}\oplus \Sigma N^{\dagger}) \t{ by \cite[Lemma 8.32]{lw}}\\
& \simeq \Hom_{\ul{\msf{CM}}(R^{\sharp\sharp})}(M, N^{\dagger}\oplus \Sigma N^{\dagger}), \t{ again by \cite[Lemma 8.32]{lw}};
\end{align*}
which shows that $(-\oplus\Sigma-)^{\dagger}$ is a right adjoint. That it is also a left adjoint follows from the fact that $\Hom_{\ul{\msf{CM}}(R^{\sharp\sharp})}((M\oplus\Sigma M)^{\dagger},A^{\dagger})\simeq \Hom_{\ul{\msf{CM}}(R)}(M,A\oplus\Sigma^{-1}A)$, and $A\oplus\Sigma^{-1}A\simeq R\otimes_{R^{\sharp\sharp}}A^{\dagger}$ by \cite[Lemma 8.32]{lw}. As $\dagger$ is an equivalence of categories, this is enough. 
\end{proof}

In fact, since $R\otimes_{R^{\sharp\sharp}}-$ is enhanced, as is $(-)^{\dagger}$, it follows that there is a biadjunction 

\[
\begin{tikzcd}[column sep=2cm]
\scr{D}_{\t{sg}}(R^{\sharp\sharp}) \arrow[r, shift left = 1ex, "R\otimes_{R^{\sharp\sharp}}"] \arrow[r, leftarrow, shift right = 1ex, swap, "\tilde{\dagger}\circ(-\oplus \Sigma)"] & \scr{D}_{\t{sg}}(R)
\end{tikzcd}
\]
on the level of small stable $\infty$-categories; for brevity write $f=\tilde{\dagger}\circ(-\oplus \Sigma)$. Since taking ind-completions is functorial, it follows that there are filtered colimit preserving functors $R\otimes_{R^{\sharp\sharp}}-\colon \msf{Ind}(\scr{D}_{\t{sg}}(R^{\sharp\sharp}))\to\msf{Ind}(\scr{D}_{\t{sg}}(R))$ and $\msf{Ind}(f)\colon \msf{Ind}(\scr{D}_{\t{sg}}(R))\to \msf{Ind}(\scr{D}_{\t{sg}}(R^{\sharp\sharp}))$.

\begin{lem}
The pair $(R\otimes_{R^{\sharp\sharp}}-,\msf{Ind}(f))$ is a biadjunction.
\end{lem}
Before proving this, recall from \cite[p. 377]{htt} that if $\scr{C}$ is a small $\infty$-category, the morphism spaces in $\msf{Ind}(\scr{C})$ are given by
\[
\Hom_{\msf{Ind}(\scr{C})}([\rlim X_{i}],[\rlim Y_{j}])=\llim_{I}\rlim_{J}\Hom_{\scr{C}}(X_{i},Y_{j}).
\]

\begin{proof}
Let $X=\rlim X_{i}$ be an object in $\msf{Ind}\,\scr{D}_{\t{sg}}(R^{\sharp\sharp})$ and $Y=\rlim Y_{j}$ an object in $\msf{Ind}\,\scr{D}_{\t{sg}}(R)$, where $X_{i}$ and $Y_{j}$ are compact objects for all $i,j$. Then there are isomorphisms
\begin{align*}
\Hom_{\msf{Ind}(\scr{D}_{\t{sg}}(R))}(R\otimes_{R^{\sharp\sharp}}X, Y) &\simeq \Hom_{\msf{Ind}(\scr{D}_{\t{sg}}(R))}(\rlim (R\otimes_{R^{\sharp\sharp}}X_{i}), \rlim Y_{j})\\
&\simeq \llim \rlim \Hom_{\msf{Ind}(\scr{D}_{\t{sg}}(R))}(R\otimes_{R^{\sharp\sharp}}X_{i}, Y_{j}) \\
&\simeq \llim \rlim \Hom_{\scr{D}_{\t{sg}}(R^{\sharp\sharp})}(X_{i},f Y_{j}) \\
&\simeq \Hom_{\msf{Ind}(\scr{D}_{\t{sg}}(R^{\sharp\sharp}))}(\rlim X_{i}, \msf{Ind}(f)(\rlim Y_{j})) \\
&\simeq \Hom_{\msf{Ind}(\scr{D}_{\t{sg}}(R^{\sharp\sharp}))}(X, \msf{Ind}(f)(Y)),
\end{align*}
which shows that $\msf{Ind}(f)$ is a right adjoint to $R\otimes_{R^{\sharp\sharp}}-$. That it is a left adjoint follows from an almost identical argument paired with \cref{functorsoncompacts}.
\end{proof}
As a consequence, both $R\otimes_{R^{\sharp\sharp}}-$ and $\msf{Ind}(f)$ are exact functors of stable $\infty$-categories.

Explicitly, we may evaluate $\msf{Ind}(f)$ on an object via
\[
\msf{Ind}(f)(\rlim Y_{i})= \rlim_{i} \tilde{\dagger}(Y_{i}\oplus \Sigma Y_{i}).
\]
As the homotopy category of $\msf{Ind}(\scr{D}_{\t{sg}}(R))$ is $\stgf{R}$, and likewise over $R^{\sharp\sharp}$, we shall let $F\colon\stgf{R}\to\stgf{R^{\sharp\sharp}}$ be the triangulated functor induced by $\msf{Ind}(f)$. This can be explicitly given on objects by letting any $X\in\stgf{R}$ be written as $X=\msf{hocolim}_{I}X_{i}$, a homotopy colimit of compacts, and thus
\[
FX=\msf{hocolim}_{I}(X_{i}\oplus\Sigma X_{i})^{\dagger}.
\]

Combining the above lemmas, we obtain the following theorem.

\begin{thm}\label{ambidextrous}
There is an infinite ladder of adjoints
\[
\begin{tikzcd}[row sep=1in, column sep=1.5cm]
\stgf{R^{\sharp\sharp}} \arrow[d, "R\otimes_{R^{\sharp\sharp}}-" description] \arrow[d, leftarrow, shift left=0.9cm, "F" description] \arrow[d, leftarrow, shift right=0.9cm, "F" description] \arrow[d, phantom, shift right=1.7cm, "\cdots"] \arrow[d, phantom, shift left=1.7cm, "\cdots"]\\
\stgf{R}.
\end{tikzcd}
\]
\end{thm}

\begin{rem}
While one obtains three adjoints from $R\otimes_{R^{\sharp\sharp}}-$ using Brown representability (a left adjoint, and a right adjoint which has its own right adjoint), it is not obvious how to obtain the infinite chain, since it is not clear that the right adjoint of $R\otimes_{R^{\sharp\sharp}}-$ preserves compact objects.
\end{rem}

\section{The case when $S$ is flat over $R$}\label{sec:flatover}
For this section, we assume that $R\to S$ is a ring map of right coherent rings, such that $S$ is a finitely presented flat module, that is finitely generated projective, over $R$. 

\begin{chunk}\label{imagegfc}
As $S$ is projective, the functor $S\otimes_{R}-$ can now be seen both as a functor $\gfc{R}\to\gfc{S}$, as previously considered, but also as a functor $\gfc{R}\to\gfc{R}$, due to the fact that $\gfc{R}$ is closed under finite direct sums and retracts, and thus tensoring with a finitely generated projective. We shall let $\msf{Im}$ denote the essential image of the functor $S\otimes_{R}-\colon\gfc{R}\to\gfc{R}$. Note that $\msf{Im}$ need not be extension closed.

As $S$ is projective, the restriction of scalars functor $\msf{res}\colon\Mod{S}\to\Mod{R}$ restricts to a functor $\msf{res}\colon\msf{GFlat}(S)\to\msf{GFlat}(R)$. This is immediate from the definition of Gorenstein flat, and the fact that restriction of scalars along a projective preserves injective modules. Moreover, the assumption that $S$ is finitely generated projective also means that $\msf{res}$ preserves cotorsion modules. This can be seen from the isomorphisms $\t{Ext}_{R}^{i}(F, \Hom_{S}(S,N))\simeq\t{Ext}_{S}^{i}(S\otimes_{R}F,N)$, which hold for each $i\geq 0$.

Consequently, we have a functor $\msf{res}\colon\gfc{S}\to\gfc{R}$. In particular, we have the adjoint pair $(S\otimes_{R}-,\msf{res})\colon\gfc{R}\to\gfc{S}$. The inclusion $\msf{Im}\to \gfc{R}$ corresponds precisely to the restriction of $\msf{res}$ to the Gorenstein flat cotorsion $S$-modules that are in the image of $S\otimes_{R}-$. 

The assumption that $S$ is projective over $R$ means $S\otimes_{R}-$ has a left adjoint, given by $S^{*}\otimes_{S}-\colon\Mod{S}\to\Mod{R}$, where $S^{*}:=\Hom_{R}(S,R)$. Indeed,
\[
\Hom_{R}(S^{*}\otimes_{S}X,Y)\simeq\Hom_{S}(X,\Hom_{R}(S^{*},Y))\simeq \Hom_{S}(X,S\otimes_{R}Y)
\] 
for any $X\in\Mod{S}$ and $Y\in\Mod{R}$, using the isomorphism $\Hom_{R}(\Hom_{R}(S,R),-)\simeq S\otimes_{R}-$ which holds as $S$ is finitely generated projective.

If $X\in\gfc{R}$, then $S^{*}\otimes_{S}(S\otimes_{R}X)\simeq S^{*}\otimes_{R}X\simeq \Hom_{R}(S,X)$. Now, viewed as a left $R$-module, $\Hom_{R}(S,X)$ is Gorenstein flat and cotorsion, and therefore there is an adjoint triple
\[
\begin{tikzcd}[column sep=1in]
\gfc{R} \arrow[r, "S\otimes_{R}-" description ] \arrow[r, leftarrow, shift left = 2ex, "S^{*}\otimes_{S}-"] \arrow[r, leftarrow, shift left = -2ex, "\msf{res}", swap] & \msf{Im}.
\end{tikzcd}
\]

We note that without further assumptions, we may not replace $\msf{Im}$ on the right hand side with $\gfc{S}$. This is because with only the assumption that $S$ is a right coherent ring, we may not deduce that $S^{*}\otimes_{S}-\colon\Mod{S}\to\Mod{R}$ preserves Gorenstein flat modules ($S^{*}$ need not have finite flat dimension over $S$) nor cotorsion modules.
\end{chunk}

\begin{chunk}\label{adjst}
Having considered the above functors on the level of Gorenstein flat cotorsion modules, let us now consider how they interact on the level of the stable categories. 

We shall let $\msf{J}$ denote the essential image of $S\otimes_{R}-\colon\stgf{R}\to\stgf{S}$. By \cref{diagram}, any object of $\msf{J}$ is an object of $\msf{Im}$.

Again, the functor $S\otimes_{R}-\colon\stgf{R}\to\msf{J}$ has a right adjoint, given by $\msf{res}$. This is clear on the level of objects, while if $f\colon S\otimes_{R}X\to S\otimes_{R}Y$ factors through a flat cotorsion $S$-module, it factors, by \cref{factorsthroughimage}, through a flat cotorsion $S$-module of the form $S\otimes_{R}F$, where $F\in\msf{FlatCot}(R)$. By the argument given in \cref{imagegfc}, $\msf{res}(S\otimes_{R}F)$ is a flat cotorsion $R$-module, and $\msf{res}(f)$ trivially factors through it.

A similar argument shows that $S^{*}\otimes_{R}-$ preserves flat cotorsion modules, and therefore also induces a functor $\msf{J}\to\stgf{R}$, which is the right adjoint to $S\otimes_{R}-$. Note that $\msf{J}$ need not be a triangulated subcategory of $\stgf{S}$, so the adjoints below are, at present, merely adjoints of additive categories.

\[
\begin{tikzcd}[column sep=1in]
\stgf{R} \arrow[r, "S\otimes_{R}-" description ] \arrow[r, leftarrow, shift left = 2ex, "S^{*}\otimes_{S}-"] \arrow[r, leftarrow, shift left = -2ex, "\msf{res}", swap] & \msf{J}.
\end{tikzcd}
\]

\end{chunk}

Let us now consider two questions that arise from the comparison between \cref{imagegfc} and \cref{adjst}: when can the adjoint triple of \cref{imagegfc} be extended to an adjoint triple between the categories of Gorenstein flat cotorsion modules? and which condition can be imposed to enables a much more explicit relationship between the adjunctions on the Frobenius level and the stable category level, particularly in relation to $\msf{Im}$ and $\msf{J}$?

The conditions to ensure a positive outcome for the first question are somewhat restrictive.

\begin{prop}\label{allofgfc}
Let $R\to S$ be a map of right coherent rings such that $S$ is finitely presented and flat over $R$. If $S^{*}=\Hom_{R}(S,R)$ is finitely presented as an $S$-module, and either $\msf{flatdim}_{S}\,S^{*}$ or $\msf{injdim}_{S}\,S^{*}$ is finite, then $S^{*}\otimes_{S}-\colon\msf{Ch}(S)\to\msf{Ch}(R)$ preserves totally acyclic complexes of flat cotorsion modules.
\end{prop}

\begin{proof}
The assumption that $S^{*}$ is finitely presented over $S$ ensures that $S^{*}\otimes_{S}-$ is a definable functor, and therefore preserves direct limits and products, and hence pure injective modules. As $S^{*}\otimes_{S}-$ has an exact right adjoint, it follows that it preserves projective modules and thus also flat modules. In particular, it preserves pure injective flat modules, which are just the flat cotorsion modules.

Consequently, if $T$ is a totally acyclic complex of flat cotorsion $S$-modules, then $(S\otimes_{S}T)_{i}\in\msf{FlatCot}(R)$ for each $i\in\mbb{Z}$. We now turn our attention to the acyclicity requirements. The assumption that $S^{*}$ has either finite flat or injective dimension ensures that $\t{Tor}_{S}^{i}(S^{*}, Z_{j}(T))=0$ for all $i>0$ and $j\in\mbb{Z}$, and thus $S^{*}\otimes_{S}T$ is an acyclic complex. 

Since the cycles of $S^{*}\otimes_{S}T$ are cotorsion, the third condition of \cref{totallyacyclicdef} is satisfies. 

Lastly, we show that $\Hom_{R}(S^{*}\otimes_{S}T,F)$ is acyclic for each $F\in\msf{FlatCot}(R)$. If $F$ is a flat cotorsion $R$-module, then there are isomorphisms

\begin{align*}
\msf{R}\Hom_{R}(S^{*} \otimes_{S}^{\msf{L}}T,F)&\simeq \msf{R}\Hom_{S}(T,\msf{R}\Hom_{R}(S^{*},F)) \\
&\simeq\msf{R}\Hom_{S}(T,\msf{R}\Hom_{R}(\msf{R}\Hom_{R}(S,R),F))\\
&\simeq\msf{R}\Hom_{S}(T,S\otimes_{R}\msf{R}\Hom_{R}(R,F))\\
&\simeq\msf{R}\Hom_{S}(T,S\otimes_{R}F)\\
&\simeq 0,
\end{align*}

where the third isomorphism follows from the fact that $S$ is projective over $R$ and the fifth follows from the fact that $S\otimes_{R}F$ is in $\msf{FlatCot}(S)$; by assumption $\msf{R}\Hom_{S}(T,X)=0$ for all $X\in\msf{FlatCot}(S)$. As a result, we see that $\Hom_{R}(S^{*}\otimes_{S}T,F)$ is acyclic, as desired, which finishes the proof.
\end{proof}

The assumption that $S^{*}$ has finite flat dimension over $S$ is not strong at all. It happens whenever we assume that either $R$ or $S$ is left Noetherian, by \cite[Proposition 8.4.14]{dcmca}, and thus particularly happens whenever $R$ and $S$ are Gorenstein rings. The condition that $S^{*}$ is finitely generated over $S$ is more restrictive, but happens in nature, as the following example shows.

\begin{ex}
Let $R$ be a Gorenstein, hence Noetherian, ring and let $P$ be a finitely generated projective $R$-module such that $P\simeq P^{*}$; this happens if, for example $P=R$, or $P=Q\oplus Q^{*}$ for any finitely generated projective module $Q$. 

Let $S=R\ltimes P$ be the trivial extension of $R$ by $P$, which, as an $R$-module, is just $R\oplus P$, and as a ring has action $(r_{1},p_{1})(r_{2},p_{2})=(r_{1}r_{2},r_{1}p_{2}+r_{2}p_{1})$. Clearly $R$ can be viewed as a subring of $S$ via $R\ltimes 0$, hence there is a natural map $R\to S$.

The module $S^{*}=\Hom_{R}(S,R)$ is just isomorphic to $S$ by construction, and thus all the conditions of \cref{allofgfc} are satisfied.

\end{ex}

\begin{cor}
Let $R\to S$ be a map of right coherent rings such that $S$ is finitely presented and flat over $R$. Suppose $S^{*}=\Hom_{R}(S,R)$ is finitely presented as an $S$-module, and either $\msf{flatdim}_{S}\,S^{*}$ or $\msf{injdim}_{S}\,S^{*}$ is finite. Then there is an adjoint triple 
\[
\begin{tikzcd}[column sep=1in]
\gfc{R} \arrow[r, "S\otimes_{R}-" description ] \arrow[r, leftarrow, shift left = 2ex, "S^{*}\otimes_{S}-"] \arrow[r, leftarrow, shift left = -2ex, "\msf{res}", swap] & \gfc{S},
\end{tikzcd}
\]
of Frobenius exact categories. This induces a corresponding adjoint triple of triangulated functors
\[
\begin{tikzcd}[column sep=1in]
\stgf{R} \arrow[r, "S\otimes_{R}-" description ] \arrow[r, leftarrow, shift left = 2ex, "S^{*}\otimes_{S}-"] \arrow[r, leftarrow, shift left = -2ex, "\msf{res}", swap] & \stgf{S}
\end{tikzcd}
\]
between the stable categories.
\end{cor}

\begin{proof}
This is just a combination of \cref{triangulatedfunctor}, \cref{allofgfc} and the discussion in \cref{imagegfc} concerning $\msf{res}$ preserving Gorenstein flat cotorsion and flat cotorsion modules.
\end{proof}

Let us now turn our attention to the second of the questions.

\begin{chunk}\label{epimorphism}
Recall that a ring map $R\to S$ is an \emph{epimorphism of rings} if the natural map $S\otimes_{R}S\to S$ is an isomorphism, see for example \cite[p. 45]{krbook} for more details.
\end{chunk}

\begin{lem}
Let $R$ and $S$ be right coherent rings and suppose that $R\to S$ is a finite flat ring epimorphism. Then $\msf{Im}$ is a Frobenius category, where the injective objects are of the form $S\otimes_{R}F$ with $F\in\msf{FlatCot}(R)$. The stable category of $\msf{Im}$ is $\msf{J}$.
\end{lem}

\begin{proof}
First we show that under these assumptions that $\msf{Im}$ is closed under extensions, and thus an exact category. If $0\to S\otimes_{R}M\to X\to S\otimes_{R}N\to 0$ is a short exact sequence in $\gfc{S}$ with $M,N\in\gfc{R}$, then, by the assumption that $S$ is projective over $R$, it follows that the same sequence is short exact in $\gfc{R}$ when viewed as a sequence of $R$-modules. As $S$ is a ring epimorphism, we have that $S\otimes_{R}X_{R}\simeq X$, and thus applying $S\otimes_{R}-$ to the restricted sequence, we see that $X\in\msf{Im}$.

For brevity, let $\msf{X}=\{S\otimes_{R}F:F\in\msf{FlatCot}(R)\}\subseteq\gfc{S}$. The assumptions of $R\to S$ means, again, that restriction of scalars preserves flat modules, and therefore, by a similar argument to above, we see that $\msf{X}$ is also an exact category. Moreover, any object of $\msf{X}$ is projective and injective in $\msf{Im}$, as it is contained within the injective objects of $\gfc{S}$.

Let us now show that $\msf{Im}$ has enough projective objects. Suppose that $S\otimes_{R}M\in\msf{Im}$, where $M\in\gfc{R}$ and let $0\to L\to F\to M\to 0$ be a short exact sequence in $\gfc{R}$ with $F\in\msf{FlatCot}(R)$. Applying $S\otimes_{R}-$ to this short exact sequence, gives a deflation $S\otimes_{R}F\to S\otimes_{R}M\to 0$, and as $S\otimes_{R}M\neq 0$ by assumption, it follows that $S\otimes_{R}F\neq 0$ and therefore, as $S\otimes_{R}F\in\msf{X}$, we see $\msf{Im}$ has enough projectives. An almost identical argument shows the existence of enough injectives.

Now suppose that $P\in\msf{Im}$ is projective. As $P\in\msf{Im}$, there is an $X\in\gfc{R}$ such that $S\otimes_{R}X=P$. Applying $S\otimes_{R}-$ to the short exact sequence $0\to \Omega_{R}^{1}(X)\to F\to X\to 0$, where $F\in\msf{FlatCot}(R)$, gives a short exact sequence $0\to S\otimes_{R}\Omega_{R}^{1}(X)\to S\otimes_{R}F\to P\to 0$ in $\msf{Im}$. By the assumption of $P$ being projective, it follows that $P$ is a summand of $S\otimes_{R}F$, and is therefore in $\msf{FlatCot}(S)$. However, $\msf{Im}$ is closed under summands, as restriction of scalars is exact and $R\to S$ is a ring epimorphism. Thus $P\in\msf{X}$. Again, an almost identical argument shows that every injective object in $\msf{Im}$ lies in $\msf{X}$. This concludes the proof of the statement that $\msf{Im}$ is Frobenius.

To see that $\ul{\msf{Im}}=\msf{J}$, it is enough to prove that $\msf{Im}$ is a full subcategory of $\stgf{S}$, as the objects of $\ul{\msf{Im}}$ and $\msf{J}$ are the same. Yet this follows immediately from \cref{factorsthroughimage}.
\end{proof}

Using the above lemma, we obtain the following.

\begin{thm}\label{epirecollement}
Let $R\to S$ be an epimorphism of right coherent rings such that $S$ is finitely generated projective over $R$. Then $S\otimes_{R}-\colon\stgf{R}\to\stgf{S}$ is a Bousfield localisation on $\stgf{R}$. In particular, there is a recollement of triangulated categories
\[
\begin{tikzcd}[column sep= 1in]
\msf{C} \arrow[r, "\mrm{inc}" description ] \arrow[r, leftarrow, shift left = 2ex, "\lambda"] \arrow[r, leftarrow, shift left = -2ex, "\rho", swap] &  \stgf{R} \arrow[r, "S\otimes_{R}-" description ] \arrow[r, leftarrow, shift left = 2ex, "S^{*}\otimes_{S}-"] \arrow[r, leftarrow, shift left = -2ex, "\msf{res}", swap] & \msf{J}=\ul{\msf{Im}}.
\end{tikzcd}
\]
\end{thm}

\begin{proof}
As $R\to S$ is a ring epimorphism, $S\otimes_{R}-$ is an endofunctor on $\stgf{R}$. Define a natural transformation $\eta\colon\t{Id}\to S\otimes_{R}$ via $X\mapsto S\otimes_{R}X$ given by $x\mapsto \1_{S}\otimes_{R}x$. The isomorphism $S\otimes_{R}S\simeq S$ shows that $\eta\circ(S\otimes_{R}-)=(S\otimes_{R}-)\circ\eta$, and that $(S\otimes_{R}-)\circ\eta\colon S\otimes_{R}-\to (S\otimes_{R}S)\otimes_{R}-$ is invertible. This shows that $S\otimes_{R}-$ is a Bousfield localisation. The existence of the recollement follows from \cite[Proposition 4.13.1]{krloc}, provided the inclusion $\msf{C}\to \stgf{R}$ admits a left adjoint. As this inclusion always preserves products, a left adjoint exists by \cref{brownrep} as Brown representability holds in the well generated triangulated category $\stgf{R}$, see the discussion in \cref{coproducts}. The fact that the image of $S\otimes_{R}-$ is $\msf{J}$ is the preceding lemma.
\end{proof}

\begin{ex}
Let us give an example of some ring maps which satisfy the conditions of the above theorem.
\begin{enumerate}
\item Let $R$ be a right coherent ring and $e\in R$ a non-trivial central idempotent. The canonical map $R\to eRe$ is a finite flat epimorphism, as $eRE=eR$ which is projective. Thus the functor $eRe\otimes_{R}-\colon\Mod{R}\to\Mod{eRe}$ induces a functor $\gfc{R}\to\gfc{eRe}$, as well as an essentially surjective functor $\stgf{R}\to\stgf{eRe}$. The kernel of this functor on triangulated categories, and thus the left hand side of the recollement in \cref{epirecollement} is $\stgf{\tilde{e}R\tilde{e}}$, where $\tilde{e}=1-e$.

\item Let $f\colon X\to Y$ be a morphism of affine schemes. Then, by \cite[17.9.1]{groth}, $f$ is an open immersion if and only if the corresponding ring homomorphism is a finitely presented flat ring epimorphism. Thus the morphisms of commutative rings that satisfy \cref{epirecollement} are precisely the ones which correspond to open immersions of affine schemes.
\end{enumerate}
\end{ex}

\subsection*{Acknowledgements}
This research was supported by the grant GA~\v{C}R 20-02760Y from the Czech Science Foundation, the project PRIMUS/23/SCI/006 from Charles University, and by the Charles University Research Center programmes UNCE/SCI/022 and UNCE/24/SCCI/022.

The author is grateful to  S. Opper, M. Prest, L. Shaul, J. Williamson and A. Zvonareva for their feedback and comments on this work throughout its various stages. I also thank the anonymous referee for their useful comments and observations, which have helped improve this work.

\bibliographystyle{abbrv}
\bibliography{references.bib}

\end{document}